%% file: CMGs-arxiv.tex
\documentclass[aos,preprint]{imsart}

\input deff3

\usepackage{enumerate}
\RequirePackage[OT1]{fontenc}
\RequirePackage{amsthm,amsmath,natbib}
\RequirePackage[colorlinks,citecolor=blue,urlcolor=blue]{hyperref}

\startlocaldefs
\endlocaldefs

\setcitestyle{numbers,square}

\begin{document}

\begin{frontmatter}
\title{Unifying Markov Properties\\ for Graphical Models}
\runtitle{Unifying Markov Properties}

\begin{aug}

\author{\fnms{Steffen} \snm{Lauritzen}
\ead[label=e1]{lauritzen@math.ku.dk}}%
\and
\author{\fnms{Kayvan} \snm{Sadeghi}%
\ead[label=e2]{k.sadeghi@statslab.cam.ac.uk}}

\address[a]{Department of Mathematical Sciences\\
University of Copenhagen\\
Universitetsparken 5\\
2100 Copenhagen\\
Denmark\\
\printead{e1}}

\address[b]{Statistical Laboratory\\
Centre for Mathematical Sciences\\
Wilberforce Road\\
Cambridge, CB3 0WA\\
United Kingdom\\
  \printead{e2}}

\runauthor{Lauritzen and Sadeghi}

\thankstext{t2}{Work of the second author was partially supported by grant $\#$FA9550-14-1-0141 from the U.S. Air Force Office of Scientific
Research (AFOSR) and the Defense Advanced Research Projects Agency (DARPA).}

\affiliation{University of Copenhagen and University of Cambridge}

\end{aug}

\begin{abstract}
Several types of graphs with different conditional independence interpretations --- also known as Markov properties --- have been proposed and used in graphical models. In this paper we unify these Markov properties by introducing a class of graphs with four types of edges --- lines, arrows, arcs, and dotted lines --- and a single separation criterion. We show that independence structures defined by this class specialize to each of the previously defined cases, when suitable subclasses of graphs are considered. In addition, we define a pairwise Markov property for the subclass of chain mixed graphs which includes chain graphs with the LWF interpretation, as well as summary graphs (and consequently ancestral graphs). We prove the equivalence of this pairwise Markov property to the global Markov property for compositional graphoid independence models.
\end{abstract}



\end{frontmatter}

\section{Introduction}\label{sec:1} 
Graphical models provide a strong and clear formalism for studying conditional independence relations that arise in different statistical contexts. Originally, graphs with a single type of edge were used; see, for example,  \cite{dar80} for undirected graphs (originating from statistical physics \cite{gibbs:02}),  and \cite{wer83,kii84} for directed acyclic graphs (originating from genetics \cite{wright:21}).

With the introduction of chain graphs \cite{lau89}, and other types of graphs with edges of several types \cite{cox93,wer94,ric02,pen14} as well as different interpretations of chain graphs \cite{and01,drt09}, a plethora of  Markov properties have emerged. These have been introduced with different motivations: chain graphs as a unification of directed and undirected graphs, the so-called AMP Markov property to describe dependence structures among regression residuals, bidirected graphs to represent structures of marginal independence, and other mixed graphs to represent selection effects and incomplete observations in causal models. Despite the similarities among these, the lack of a general theory  as well as the use of different definitions and notation  has undermined the original conceptual simplicity of graphical models. This motivates a unification of the corresponding Markov properties. In \cite{sadl14}, we attempted this for different types of mixed graphs, but failed to include chain graph Markov properties. Here we follow an analogous approach using a single separation criterion, but using four types of edges: line, arrow, arc, and dotted line. To the best of our knowledge, this unifies most graphical independence models previously discussed in the literature. One exception is Drton's \cite{drt09} type III chain graph Markov property which has several unfortunate properties and so far has not played any specific role; we have chosen  to avoid introducing a fifth type of edge to accommodate this property; another exception is the \emph{reciprocal} graphs of Koster \citep{koster:96}, which allow feedback cycles; 
other exceptions use graphs to describe conditional independence in dynamical systems \citep{eichler:07,didelez:08} which we do not discuss here.
Our unification includes summary graphs --- which include ancestral graphs as well as chain graphs with the multivariate regression Markov property \cite{cox93} --- chain graphs with the LWF Markov property \cite{lau89,fry90}, and chain graphs with the AMP Markov property \cite{and01}.

%

In addition to the unification of the (global) Markov property, we provide a unified pairwise Markov property. However, it seems technically complex to include the pairwise Markov property for chain graphs with the AMP interpretation and hence we only discuss this for the subclass of graphs with three types of edges where cycles of specific types are absent.  Such graphs were called chain mixed graphs (CMGs) in \cite{sad16} and its corresponding independence model unifies those of summary graphs (and ancestral graphs) as well as chain graphs with the LWF Markov property.  For CMGs, we first discuss the notion of maximality and show that every missing edge in a maximal CMG corresponds to an independence statement, thus forming a potential base for specifying pairwise Markov properties.  %
For CMGs we prove the equivalence of pairwise and global Markov properties for abstract independence models which are compositional graphoids. 

The structure of the paper is as follows: In the next section, we define graphs with four types of edges and provide basic graph theoretical definitions. In Section \ref{sec:3}, we discuss general independence models and compositional graphoids, provide a single separation criterion for such graphs, and show that the induced independence models are compositional graphoids. Further we demonstrate how the various independence models discussed in the literature are represented within this unification. In Section \ref{sec:4}, we define the notion of maximal graphs,  provide conditions under which a CMG is maximal, and show that any CMG can be modified to become maximal without changing its independence model. 
In Section \ref{sec:4n}, we  provide a pairwise Markov property for CMGs, and prove that for compositional graphoids, the pairwise Markov property is equivalent to the global Markov property. 
Finally, we conclude the paper with a discussion in Section \ref{sec:5n}.
\section{Graph terminology}\label{sec:2}
\subsection{Graphs} A \emph{graph} $G$ is a triple consisting of a \emph{node} set or
\emph{vertex} set $V$, an \emph{edge} set $E$, and a relation that with
each edge associates two nodes (not necessarily distinct), called
its \emph{endpoints}. When nodes $i$ and $j$ are the endpoints of an
edge, these are
\emph{adjacent} and we write $i\sim j$. We say the edge is \emph{between} its two
endpoints. We usually refer to a graph as an ordered
pair $G=(V,E)$. Graphs $G_1=(V_1,E_1)$ and $G_2=(V_2,E_2)$ are called \emph{equal} if $(V_1,E_1)=(V_2,E_2)$. In this case we write $G_1=G_2$.

The graphs that we use 
are \emph{labeled graphs}, i.e.\ every node is considered a different object. Hence, for example, the graph $i\ful j\ful k$ is not equal to the graph $j\ful i\ful k$.

In addition, in this paper, we use graphs with four types of edges denoted by
\emph{arrows}, \emph{arcs} (solid lines with two-headed arrows), \emph{lines} (solid lines), and \emph{dotted lines}; as will be seen in Section \ref{sec:3}, we shall use dotted lines to represent chain graphs with the AMP Markov property. Henceforth, by `graph', we mean a graph with these four possible types of edges. We do not distinguish between $i\ful j$ and $j\ful i$, between $i\arc j$
and $j \arc i$, or between $i\ddash j$ and $j\ddash i$, but we do
distinguish between $j\fra i$ and $i\fra j$.

A \emph{loop} is an edge with
endpoints being identical. In this paper, we are only considering graphs that do not contain loops. \emph{Multiple edges} are edges sharing the
same pair of endpoints. A \emph{simple graph} has neither loops nor multiple edges.
Graphs we are considering in this paper may generally contain multiple edges, even of the same type. However we shall emphasize for all purposes in the present paper, multiple edges of the same type are redundant and hence at most one edge of every type is necessary to represent the objects we discuss.

We say that $i$ is a
\emph{neighbor} of $j$ if these are endpoints of a line; if there is an arrow from $i$ to $j$, $i$ is a \emph{parent} of $j$ and $j$ is a \emph{child} of $i$. We also say that $i$ is a \emph{spouse} of $j$ if these are endpoints of an arc, and $i$ is a \emph{partner} of $j$ if they are endpoints of a dotted line. We use the notations $\nei(j)$, $\pa(j)$, $\spo(j)$, and $\partn(j)$ for the set of all neighbours, parents, spouses, and partners of $j$ respectively. More generally, for a set of nodes $A$ we let $\nei(A)= \cup_{j\in A} \nei(j)\setminus A$ and similarly for $\pa(A)$, $\spo(A)$, and $\partn(A)$.


A \emph{subgraph} of a graph $G_1$ is graph $G_2$ such that $V(G_2)\subseteq V(G_1)$ and each edge present in $G_2$ also occurs in $G_1$ and has the same
type there. An \emph{induced subgraph} by a subset $A$ of the node set is a subgraph that contains all and only nodes in $A$ and all edges between two
nodes in $A$.

A \emph{walk} $\omega$ is a list $\omega=\langle i_0,e_1,i_1,\dots,e_n,i_n\rangle$ of nodes and edges such that for $1\leq m\leq n$, the edge $e_m$ has endpoints $i_{m-1}$ and $i_m$. We allow a walk to consist of a single node $i_0=i_n$. If the graph is simple then a walk can be determined uniquely by a sequence of nodes. Also, a non-trivial walk is always determined by its edges, so we may write $\omega=\langle e_1,\dots,e_n\rangle$ without ambiguity. Throughout this paper, however, we often use only node sequences to describe walks even in graphs with multiple edges, when it is apparent from the context or the type of the walk which edges are involved.
The first and the last nodes  of a walk are its \emph{endpoints}. All other nodes are  \emph{inner nodes} of the walk. We say a walk  is \emph{between} its endpoints. A \emph{cycle} is a walk with at least
two edges and no repeated node except $i_0=i_n$. A \emph{path} is a walk with no repeated node.

A \emph{subwalk} of a walk $\omega=\langle i_0,e_1,i_1,\dots,e_n,i_n\rangle$ is a walk that is a subsequence $\langle i_r,e_{r+1},i_{r+1},\dots,e_p,i_p\rangle$ of $\omega$ between two occurrences of nodes ($i_r,i_p$, $0\leq r\leq p\leq n$). If a subwalk forms a path then it is also a \emph{subpath} of $\omega$.


In this paper we need different types of walks as defined below. Consider a walk $\omega=\langle i=i_0,i_1,\dots,i_n=j\rangle$. We say that
\begin{itemize}
  \item $\omega$ is \emph{undirected} if it only consists of solid lines;
  \item $\omega$ is \emph{directed} from $i$ to $j$ if all  edges $i_qi_{q+1}$, $0\leq q\leq n-1$, are arrows pointing from $i_q$ to $i_{q+1}$;
  \item $\omega$ is \emph{semi-directed} from $i$ to $j$ if it has at least one arrow, no arcs, and every arrow $i_qi_{q+1}$ is pointing from $i_q$ to $i_{q+1}$;
  \item $\omega$ is \emph{anterior} from $i$ to $j$ if it is semi-directed from $i$ to $j$ or  if it is composed of lines and dotted lines.
\end{itemize}
Thus a directed walk is also semi-directed and a semi-directed walk is also an anterior walk.
If there is a directed walk from $i$ to $j$ ($j\neq i$) then $i$ is an \emph{ancestor} of $j$. We denote the set of ancestors of $j$ by $\an(j)$. If there is an anterior walk from $i$ to $j$ ($j\neq i$) then we also say that $i$ is \emph{anterior} of $j$. We use
the notation $\ant(j)$ for the set of all anteriors of $j$. For a set $A$, we define $\ant(A)=\bigcup_{i\in A}\ant(i)\setminus A$. We also use the notations $\An(A)$ and $\Ant(A)$ for the set of \emph{reflexive} ancestors and anteriors of $A$ so that $\An(A)=A\cup\an(A)$ and $\Ant(A)=A\cup\ant(A)$. In addition, we define a  set $A$ to be \emph{anterior} if $\ant(i)\subseteq A$ for all $i\in A$; in other words, $A$ is anterior if $\ant(A)=\varnothing$.

In fact, we are only interested in these walks when we discuss graphs without dotted lines. For example, consider the following walk (path) in such a graph:
$$i\ful j\ful k\fra l\fra m\ful n\fra o\arc p.$$
Here it holds that there is an undirected walk between $i$ and $k$ and hence $i\in\ant(k)$, but there is no semi-directed walk from $i$ to $k$. In addition, we have that $k\in\an(m)$ and $i\in\ant(o)$, while there is a semi-directed walk from $i$ to $o$. There is also no anterior walk from $i$ to $p$.


Notice that, unlike most places in the literature (e.g.\  \cite{ric02}), we use walks instead of paths to define ancestors and anteriors.  Using walks instead of paths is immaterial for this purpose as the following lemma shows.

\begin{lemma}\label{lem:1130n}
There is a directed or anterior walk from $i$ to $j$ if and only if there is a directed or anterior path from $i$ to $j$ respectively.\end{lemma}
\begin{proof} If there is a path, there is a walk as a path is also a walk. Conversely, assume there is a directed or anterior walk from $i$ to $j$. If $i=j$ then we are done by definition.
Otherwise, start from $i$ and move on the walk towards $j$. Consider the first place where a node $k$ is repeated on the walk. The walk from $k$ to $k$ forms a cycle. If we remove this cycle from the walk, the resulting walk remains directed; similarly, the walk resulting from an anterior walk remains anterior. Successively removing all cycles along the walk in this way implies the result.
\end{proof}

A \emph{section} $\rho$ of a walk  is a maximal subwalk consisting only of solid lines, meaning that there is no other subwalk that only consists of solid lines and includes $\rho$. A walk decomposes uniquely into sections; 
sections may also be single nodes.  
The section is an \emph{inner section} on the walk if all nodes on the section are inner nodes on the walk and an \emph{endpoint section} if it contains an endpoint of the walk.
A section $\rho$ on a walk $\omega$ is called a  \emph{collider section} if one of the four following walks is a subwalk of $\omega$: $u\fra\rho\fla\,v$, $u\arc\rho\fla\,v$, $u\arc\rho\arc\,v$, $u\fra\rho\ddash\,v$, and $u\arc\rho\ddash\,v$, i.e., a section $\rho$ is a collider if two arrowheads meet at $\rho$ or an arrowhead meets a dotted line. All other sections on $\omega$ are called \emph{non-collider} sections; these are sections that are an endpoint of $\omega$ or the following sections: $u\fla\rho\fra\,v$, $u\fla\rho\arc\,v$, $u\fla\rho\ddash\,v$, $u\fra\rho\fra\,v$, and $u\ddash\rho\ddash\,v$. We may speak of collider or non-collider sections (or nodes) without mentioning the relevant walk when this is apparent from the context. Notice that a section may be a collider on one part of the walk and a non-collider on another. For example, in Fig.~\ref{fig:acyclic graphex}(a), the section $\langle h,q\rangle$ is a collider on the walk $\langle l,h,q,p\rangle$. It is also a collider on $\langle k,q,h,p\rangle$ via the edge $h\arc p$, but a non-collider on $\langle k,q,h,p\rangle$ via the edge $h\fra p$. Notice also that $\langle k\rangle$ is a non-collider on $\langle j,k,q\rangle$.

A \emph{tripath} is a path with three distinct nodes. Note that \cite{sad13} used the term V-configuration for such a path. 
If the inner node on a tripath is a collider we shall also say that the tripath itself is a collider or non-collider. 

\subsection{Subclasses of graphs}\label{sec:2.1}
Most graphs discussed in the literature are subclasses of the graphs considered here. In addition, the global Markov property defined in the next section specializes to the independence structures previously discussed. Exceptions include \emph{MC graphs} \cite{kos02} and \emph{ribbonless graphs} \cite{sad13}. However, any independence structure represented by an MC graph or a ribbonless graph can also  be represented by a summary graph or an ancestral graph \citep{sadl14}, which are also covered in this paper.

\begin{figure}
\centering
\begin{tabular}{cc}
\scalebox{0.15}{\includegraphics{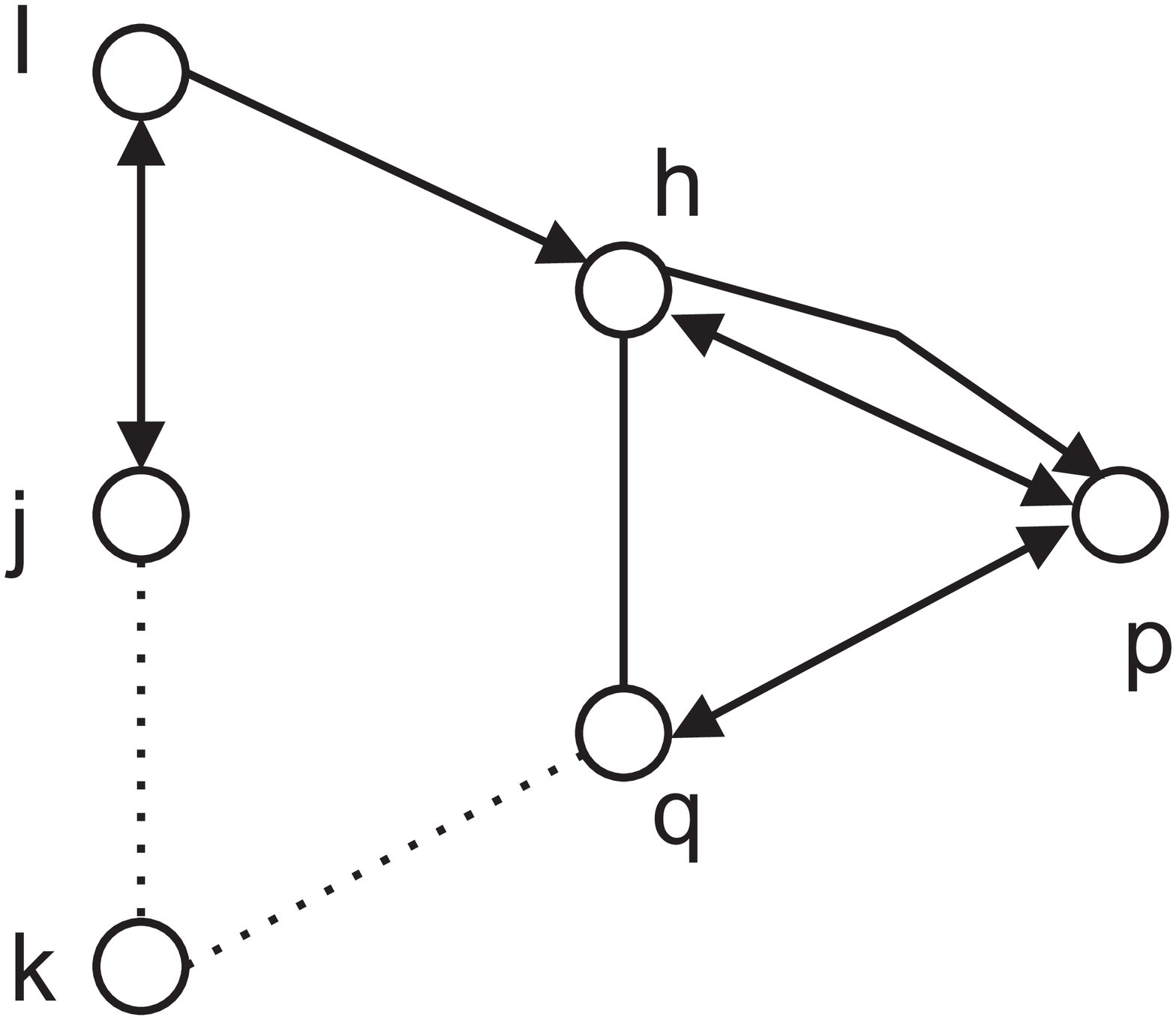}}\nn\nn &
\nn\nn\scalebox{0.15}{\includegraphics{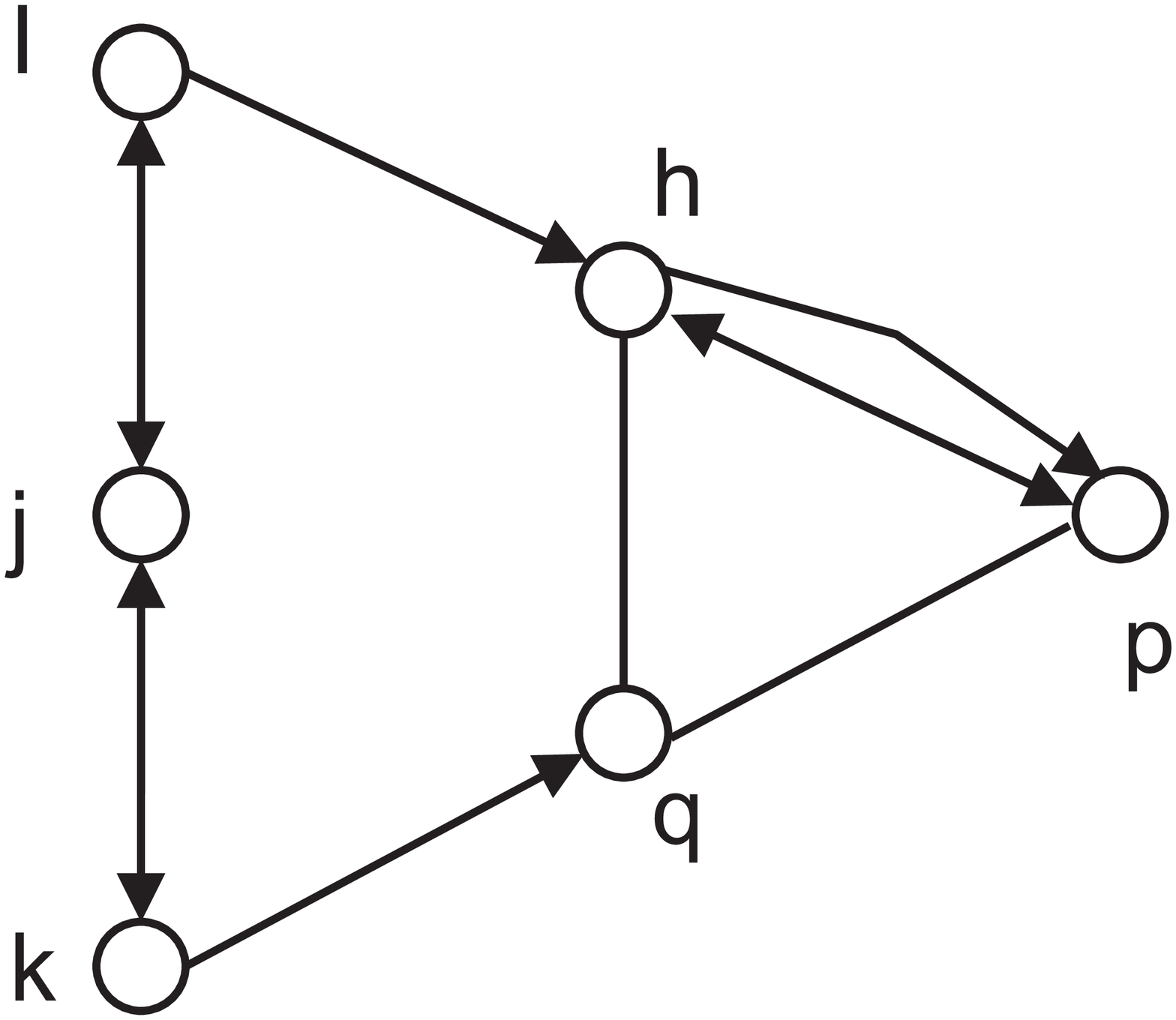}}\\
(a) & (b)
\end{tabular}
  \caption[]{\small{(a) A graph with four types of edges. (b) A graph that is not a CMG as $\langle h\fra p\ful q\ful h\rangle$ is a semi-directed cycle.}}
     \label{fig:acyclic graphex}
\end{figure}

Although we do not set any constraints on the class of graphs with four types of edges for the purpose of defining a global Markov property in Section \ref{sec:3}, the most general class of graphs for which we explicitly define a pairwise Markov property in Section \ref{sec:4n} is the class of \emph{chain mixed graphs} (CMGs) \cite{sad16}. CMGs are graphs without dotted lines and semi-directed cycles, hence \emph{reciprocal graphs} as in \cite{koster:96} are not CMGs. CMGs may have multiple edges of all types except a combination of arrows and lines or arrows in opposite directions as such combinations would constitute semi-directed cycles. The graph in Fig.\ \ref{fig:acyclic graphex}(a) is an example of graph with four types of edges, and the graph in Fig.\ \ref{fig:acyclic graphex}(b) is not a CMG because of the  semi-directed cycle $\langle h\fra p\ful q\ful h\rangle$.

It is helpful to classify subclasses of graphs into three categories: \emph{basic graphs, chain graphs, and mixed graphs,} as briefly described below.

\paragraph{Basic graphs} These are graphs that only contain one type of edge; they include \emph{undirected graphs} (UGs), containing only lines;  \emph{bidirected graphs} (BGs), containing  only bidirected edges; \emph{dotted line graphs} (DGs), containing only dotted lines; and \emph{directed acyclic graphs} (DAGs), containing only arrows without any directed cycle. Clearly, a graph without arrows has no semi-directed cycles, and a semi-directed cycle in a graph with only arrows is a directed cycle. Note that \citep{cox93,kau96,wer98,drt08}  use the terms \emph{concentration graphs} and  \emph{covariance graphs} for UGs and BGs, referring to their independence interpretation associated with covariance and concentration matrices for Gaussian graphical models. DGs have not been studied specifically; as we shall see, any independence structure associated with a DG is Markov equivalent to the corresponding UG, where dotted lines are replaced by lines. DAGs have in particular been useful  to describe causal Markov relations; see for example \citep{kii84,pea88,lauritzen:spiegelhalter:88,gei90,Spi00}. 

\paragraph{Chain graphs} A \emph{chain graph} (CG) is a graph with the two following properties: 1) if we remove all arrows, all connected components of the resulting graph --- called \emph{chain components} --- contain one type of edge only; 2) if we replace every chain component by a node then the resulting graph is a DAG.  DAGs, UGs, DGs, and BGs are all instances of chain graphs. For a DAG, all chain components are singletons, and for a chain graph without arrows, the chain components are simply the connected components of the graph.

\begin{figure}[htb]
\centering
\begin{tabular}{ccc}
\scalebox{0.15}{\includegraphics{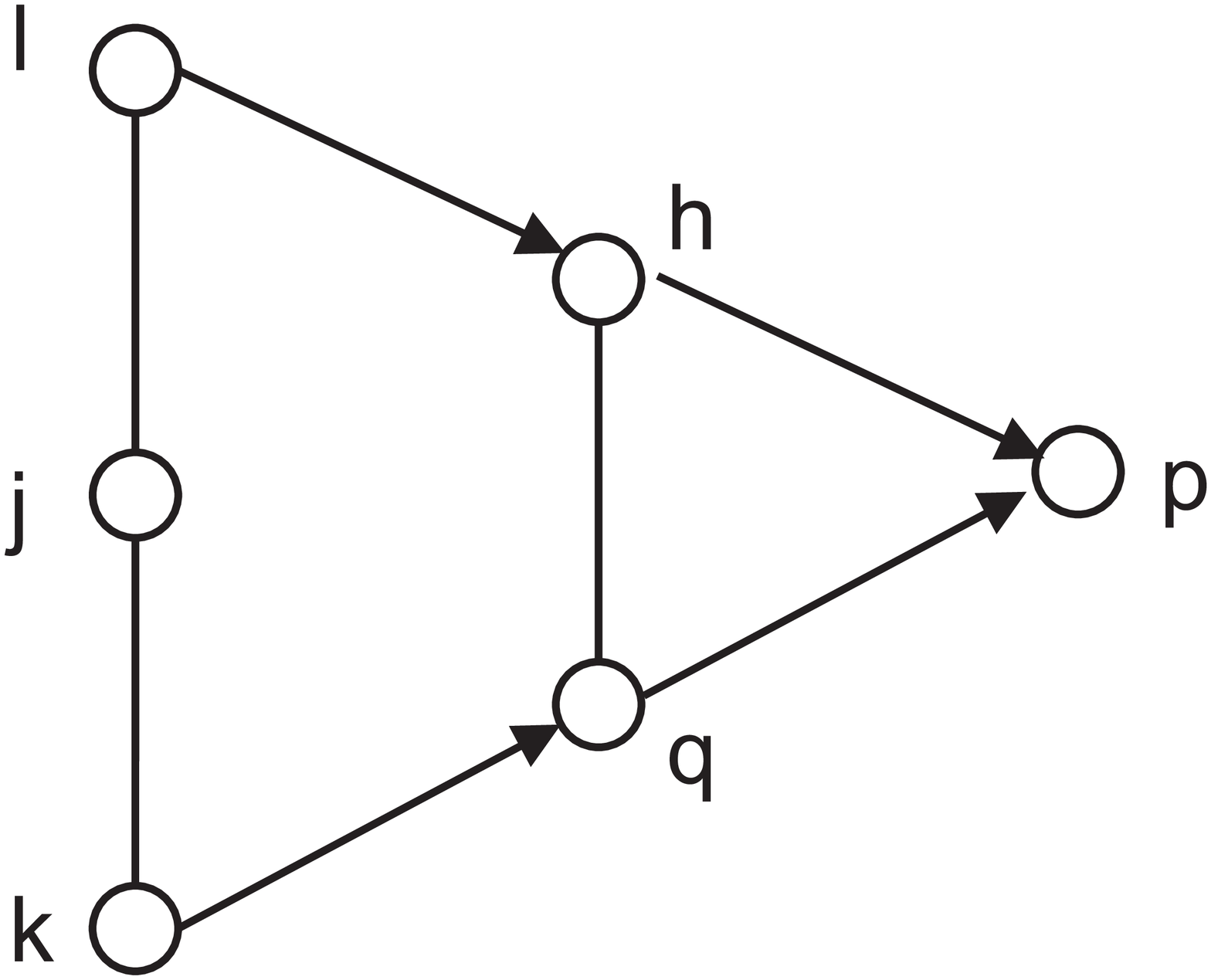}}\nn\nn &
\nn\nn\scalebox{0.15}{\includegraphics{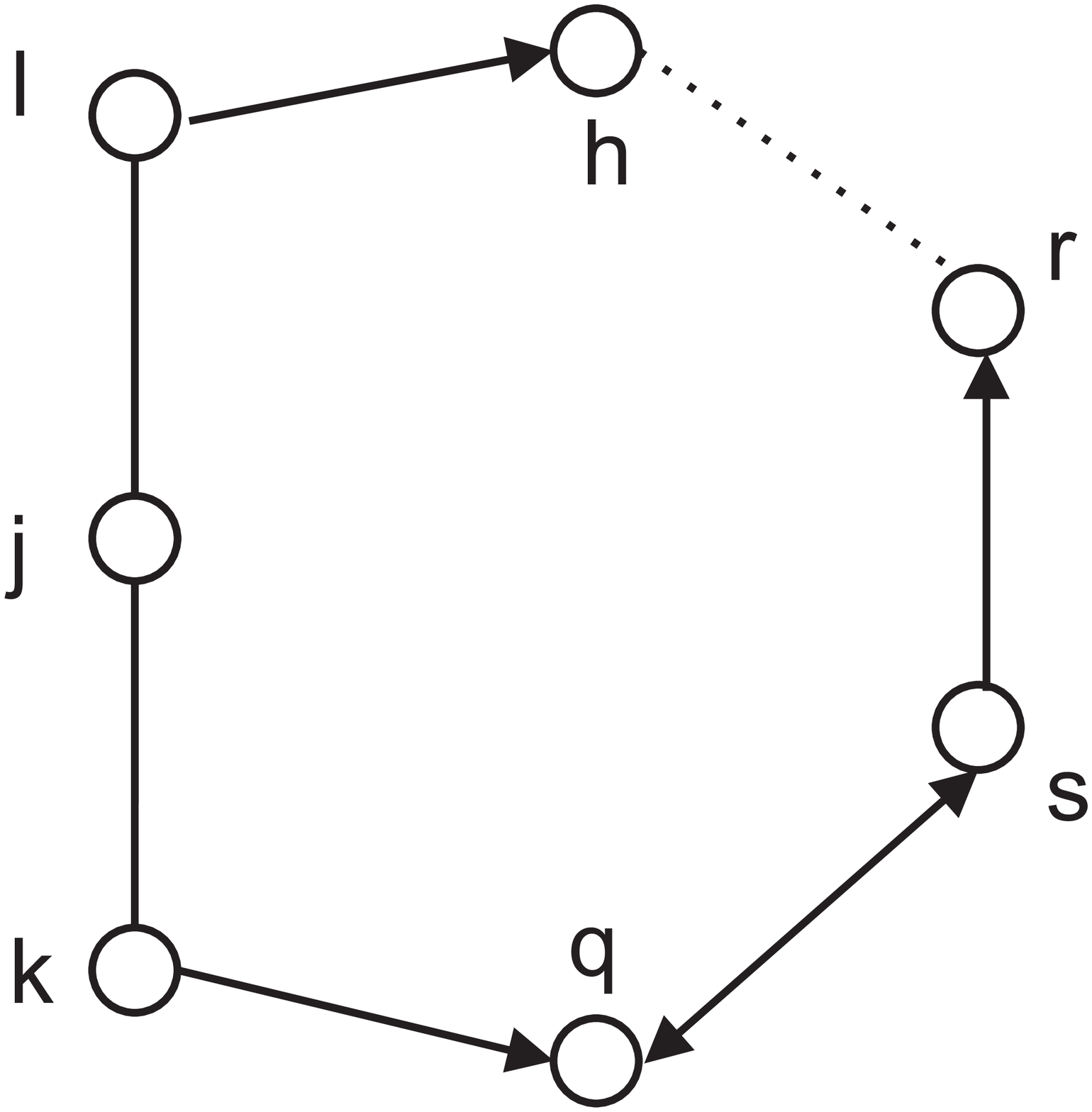}}\nn\nn&
\nn\nn\scalebox{0.15}{\includegraphics{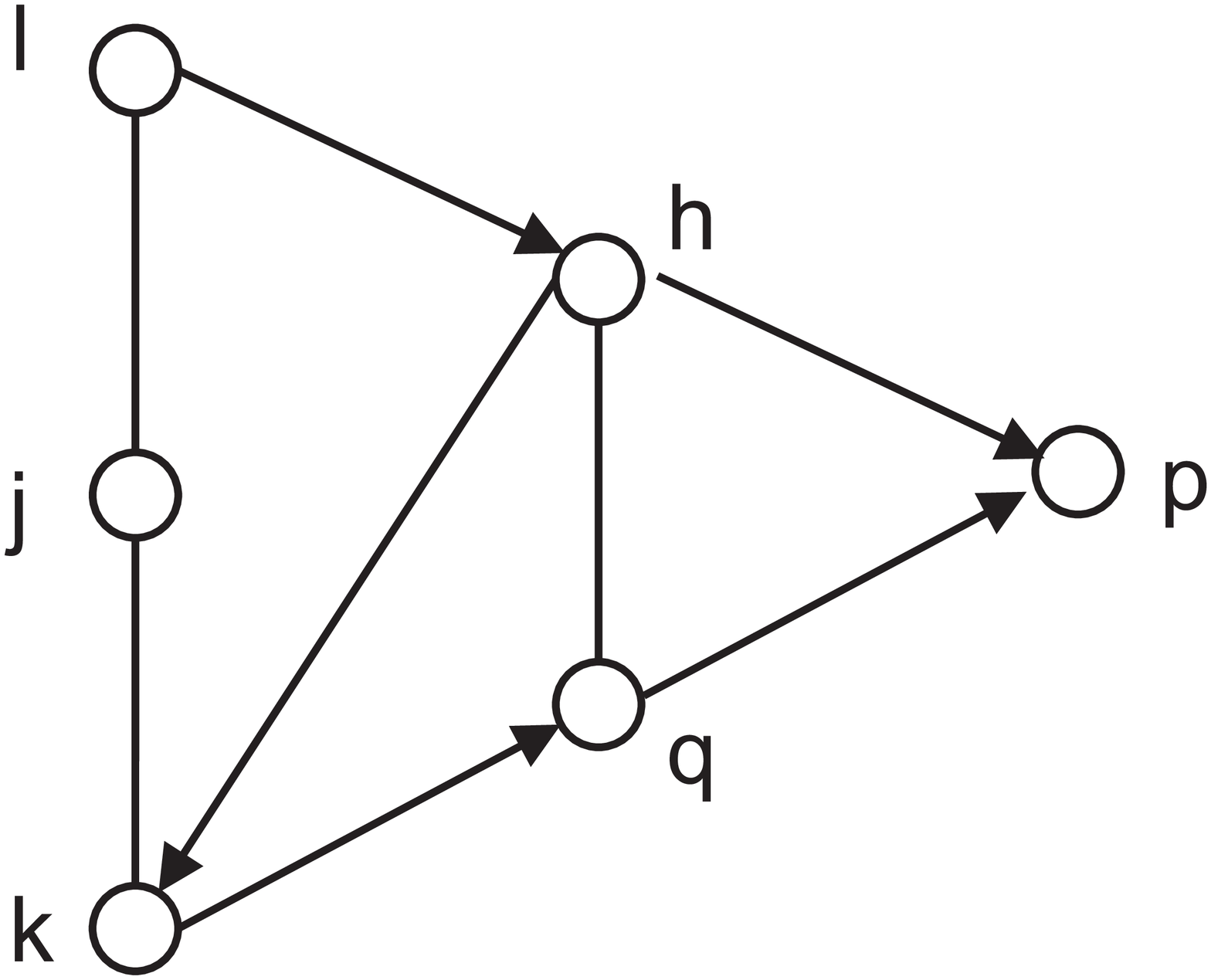}}\\
(a) & (b) & (c)
\end{tabular}
  \caption[]{\small{(a) An undirected chain graph. (b) A chain graph with chain components of different types. (c) A graph that is not a chain graph as $\langle h,k,q,h\rangle$ is semi-directed cycle in this graph.}}
     \label{fig:chainex}
\end{figure}

If all chain components contain lines, the chain graph is an \emph{undirected chain graph} (UCG) (here associated with the LWF Markov property); if all contain arcs,  it is a \emph{bidirected chain graph} (BCG) (here associated with the multivariate regression chain graph Markov property); and if all contain dotted lines, it is a \emph{dotted line chain graph} (DCG) (here associated with the AMP Markov property). 
%
%
%
For example, in Fig.\ \ref{fig:chainex}(a) the graph is a chain graph with chain components $\tau_1=\{l,j,k\}$, $\tau_2=\{h,q\}$, and $\tau_3=\{p\}$, but in Fig.\ \ref{fig:chainex}(c) the graph is not a chain graph because of the  semi-directed cycle $\langle h,k,q,h\rangle$.

\emph{Regression graphs} \cite{wers11} are chain graphs consisting of lines and arcs (although dashed undirected edges have previously been used instead of arcs in the literature), where there is no arrowhead pointing to nodes that are endpoints of lines.
\paragraph{Mixed graphs}
Marginalization and conditioning in DCGs (studied in \cite{pen14}) lead to \emph{marginal AMP graphs} (MAMPs); in our formulation, where we use dotted lines in place of full lines, MAMPs are  graphs without solid lines that satisfy three additional conditions:
\begin{enumerate}
\item $G$ has no \emph{quasi-directed cycles} in the sense it has no walk $\langle i = i_0 , i_1 , \dots , i_n = i\rangle$ containing at least one arrow and every arrow $i_qi_{q+1}$ is pointing from $i_q$ to $i_{q+1}$;
\item $G$ has no cycles composed of dotted lines and one arc;
\item If $i\ddash j\ddash k$ and $j\arc l$ for some $l$, then $i\ddash k$.
\end{enumerate}
Graphs discussed here also contain different types of \emph{mixed graphs}, a term previously used to denote graphs with lines, arrows, and arcs. These were introduced to describe independence structures obtained by marginalization and conditioning in DAG independence models; see for example \cite{sad13} for a general discussion of this issue. Examples are \emph{summary graphs} (SGs) \citep{wer11}, \emph{ancestral graphs} (AGs) \citep{ric02} and \emph{acyclic directed mixed graphs} (ADMGs) \citep{spi97,ric03}.
Summary graphs are CMGs that have no arrowhead pointing to nodes that are endpoints of lines. Ancestral graphs satisfy in addition that there are no arcs with one endpoint being an ancestor of the other endpoint. Note that in many papers about summary graphs, dashed undirected edges have  been used in place of bidirected edges.

ADMGs are summary graphs without lines. \emph{Alternative ADMGs} (AADMGs) were defined in \cite{pen16}, where arcs in ADMGs were replaced by dotted lines with our notation, although lines were used in the original definition.

CMGs are also mixed graphs, and originally defined in \cite{sad16} in order to describe independence structures obtained by marginalization and conditioning in  chain graph independence models.  \emph{Anterial graphs} (AnGs) were also defined in \cite{sad16} for the same purpose, and they are CMGs in which an endpoint of an arc cannot be an anterior of the other endpoint.

\begin{figure}[htb]
\centering
\scalebox{0.35}{\includegraphics{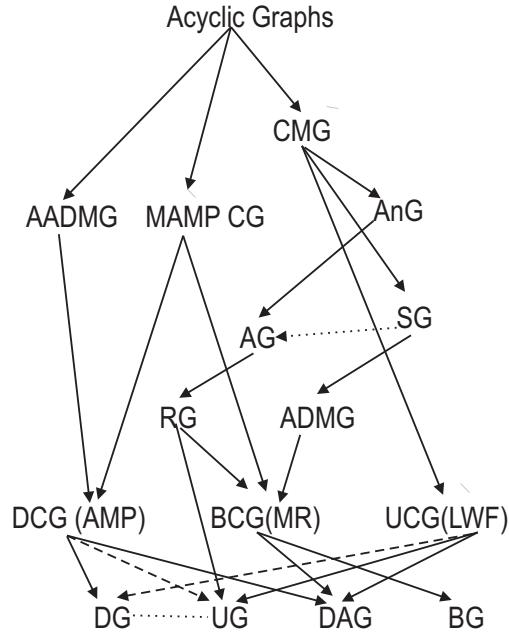}}
  \caption[]{\small{The hierarchy of graphs with four types of edges and their independence models.}}
     \label{fig:hierarchy}
\end{figure}

The diagram in Fig.~\ref{fig:hierarchy} illustrates the hierarchy of subclasses of graphs with four types of edges. 
Below we shall provide a unified separation criterion for all graphs with four types of edges and thus the associated independence models share the same hierarchy. The diagram is to be read transitively in the sense that, for example, BGs are also AGs, since the class of BGs form a subclass of BCGs, which again form a subclass of AGs; thus we omit the corresponding arrow from AG to BG.

The dashed arrow from DCG to UG indicates that although UGs are not DCGs, their associated independence models contain all independence models given by UGs and similarly for the dashed arrow from UCG to DG. The dotted arrow from SG to AG indicates that although AG is a subclass of SG, their associated independence models are the same. The dotted link between UG and DG indicates that the associated independence models are the same. These facts will be demonstrated in the next section.

\section{Graphical independence models}\label{sec:3}Graphs are used to encode independence structures for graphical models; in this section
we shall demonstrate how this can be done.
\subsection{Independence models and compositional graphoids}
An \emph{independence model} $\mathcal{J}$ over a finite set $V$ is a set of triples $\langle A,B\cd C\rangle$ (called \emph{independence statements}), where $A$, $B$, and $C$ are disjoint subsets of $V$; $C$
may be empty, but $\langle \varnothing,B\cd C\rangle$ and $\langle A,\varnothing\cd C\rangle$ are always included in $\mathcal{J}$. The independence statement $\langle A,B\cd C\rangle$ is read as ``$A$ is independent of $B$ given $C$''. Independence models may have a  probabilistic interpretation---see Section~\ref{sec:prob} for details---but this need not necessarily be the case. Similarly, not all independence models can be easily represented by graphs. For further discussion on general independence models, see \cite{stu05}.

An independence model $\mathcal{J}$ over a set $V$ is a \emph{semi-graphoid} if it satisfies the four following properties for disjoint subsets $A$, $B$, $C$, and $D$ of $V$:
 \begin{enumerate}[(S1)]
    \item $\langle A,B\cd C\rangle\in \mathcal{J}$ if and only if $\langle B,A\cd C\rangle\in \mathcal{J}$ (\emph{symmetry});
    \item if $\langle A,B\cup D\cd C\rangle\in \mathcal{J}$ then $\langle A,B\cd C\rangle\in \mathcal{J}$ and $\langle A,D\cd C\rangle\in \mathcal{J}$ (\emph{decomposition});
    \item if $\langle A,B\cup D\cd C\rangle\in \mathcal{J}$ then $\langle A,B\cd C\cup D\rangle\in \mathcal{J}$ and $\langle A,D\cd C\cup B\rangle\in \mathcal{J}$ (\emph{weak union});
    \item $\langle A,B\cd C\cup D\rangle\in \mathcal{J}$ and $\langle A,D\cd C\rangle\in \mathcal{J}$
    if and only if $\langle A,B\cup D\cd C\rangle\in \mathcal{J}$ (\emph{contraction}).
 \end{enumerate}
A semi-graphoid for which the reverse implication of the weak union property holds is said to be a \emph{graphoid}; that is it also satisfies
\begin{enumerate}[(S5)]
	\item if $\langle A,B\cd C\cup D\rangle\in \mathcal{J}$ and $\langle A,D\cd C\cup B\rangle\in \mathcal{J}$ then $\langle A,B\cup D\cd C\rangle\in \mathcal{J}$ (\emph{intersection}).
\end{enumerate}
Furthermore, a graphoid or semi-graphoid for which the reverse implication of the decomposition property holds is said to be \emph{compositional}, that is it also satisfies
\begin{enumerate}[(S6)]
	\item if $\langle A,B\cd C\rangle\in \mathcal{J}$ and $\langle A,D\cd C\rangle\in \mathcal{J}$ then $\langle A,B\cup D\cd C\rangle\in \mathcal{J}$ (\emph{composition}).
\end{enumerate}

\subsection{Independence models induced by graphs}
The notion of separation is fundamental for using graphs to represent models of independence. For a simple, undirected graph, separation has a direct intuitive meaning, so that a set $A$ of nodes is separated from a set $B$ by a set $C$ if all walks from $A$ to $B$ intersect $C$. Notice that simple separation in an undirected graph will trivially satisfy all of the properties (S1)--(S6) above, and hence compositional graphoids are  abstractions of independence models given by separation in undirected graphs. For more general graphs, separation may be more subtle, to be elaborated below.

We say that a walk $\omega$ in a graph is \emph{connecting} given $C$ if all collider sections of $\omega$ intersect $C$ and all non-collider sections are disjoint from $C$. For pairwise disjoint subsets $\langle A,B,C\rangle$, we say that $A$ and $B$ are \emph{separated} by $C$ if there are no connecting walks between $A$ and $B$ given $C$, and we use the notation $A\dse B\cd C$. The set $C$ is called an \emph{$(A,B)$-separator}.

The notion of separation above is a generalization of the \emph{$c$-separation} for UCGs as defined in \cite{stu98,stub98}.
The idea of using walks to simplify the separation theory was proposed by \cite{kos02}, who showed that, for DAGs, this notion of separation  was identical to  $d$-separation \cite{pea88}.

\begin{figure}[htb]
\centering
\scalebox{0.25}{\includegraphics{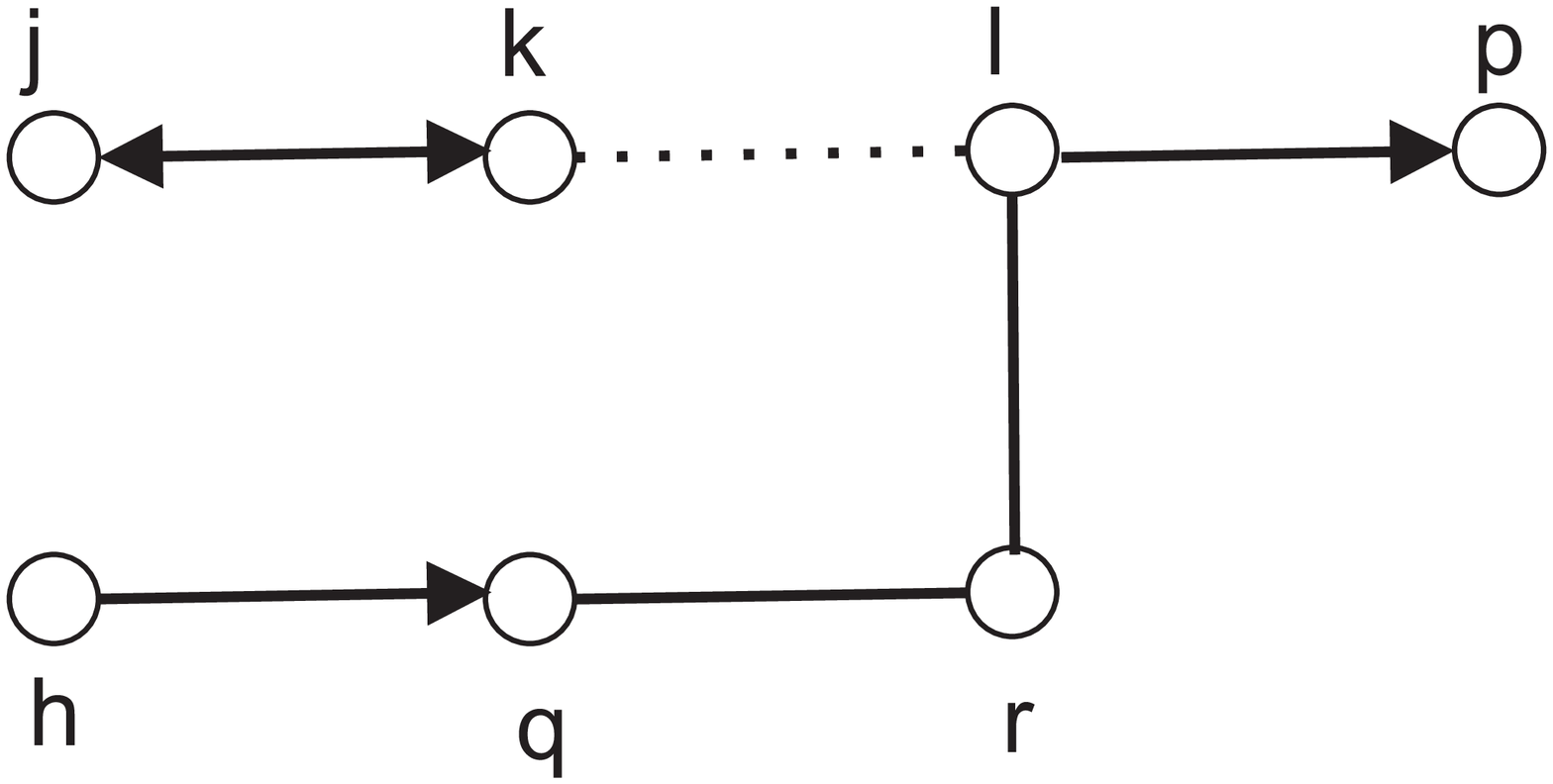}}
  \caption[]{\small{Illustration of separation in a graph $G$; it holds that $j\notdse h\cd \{k,l\}$ and  $j\notdse h\cd \{k,p\}$, but $j\dse h\cd l$ and $j\dse h\cd k$.}}
     \label{fig:sepex}
\end{figure}

For example, in the graph of Fig.\ \ref{fig:sepex}, $j\dse h\cd \{k,l\}$ and $j\dse h\cd \{k,p\}$ do not hold. The former can be seen by looking at the connecting walk $\langle j,k,l,r,q,h\rangle$, where the only node $k$ and the node $l$ of the collider sections $\langle k\rangle$ and $\langle l,r,q\rangle$ are in the potential separator set $\{k,l\}$. The latter can be seen by looking at the connecting walk $\langle j,k,l,p,l,r,q,h\rangle$, where the non-collider sections $\langle l\rangle$ and $\langle l,r,q\rangle$ are outside $\{k,p\}$, but collider sections (nodes) $\langle k\rangle$ and $\langle p\rangle$ are inside $\{k,p\}.$ However, for example, $j\dse h\cd l$ and $j\dse h\cd k$ since, in the former
case, collider section $\langle k\rangle$ is blocking all the walks and, in the latter case, one of the
collider sections $\langle l,r,q\rangle$ or $\langle p\rangle$ is blocking any walk.

A graph $G$ induces an independence model $\mathcal{J}(G)$ by separation, letting $\langle A,B\cd C\rangle\in \mathcal{J}(G)\iff A\dse B\cd C$.
It turns out that any independence model defined in this way shares the six fundamental properties of undirected graph separation. More precisely we have the following:
\begin{theorem}\label{thm:110}
For any graph $G$, the independence model $\mathcal{J}(G)$ is a compositional graphoid.
\end{theorem}
\begin{proof}
Let $G=(V,E)$, and consider disjoint subsets $A$, $B$, $C$, and $D$ of $V$. We  verify each of the six properties separately.

1) \emph{Symmetry:} If $A\dse B\cd C$ then $B\dse A\cd C$: If there is no connecting walk between $A$ and $B$ given $C$ then there is no connecting walk between $B$ and $A$ given $C$.

2) \emph{Decomposition:} If $A\dse (B\cup D)\cd C$  then $A\dse D\cd C$: If there is no connecting walk between $A$ and $B\cup D$ given $C$ then there is \emph{a forteriori} no connecting walk between $A$ and $D\subseteq (B\cup D)$ given $C$.

3) \emph{Weak union:} If $A\dse (B\cup D)\cd C$  then $A\dse B\cd (C\cup D)$: Using decomposition 2) yields $A\dse D\cd C$ and $A\dse B\cd C$. Suppose, for contradiction, that there exists a connecting walk $\omega$ between $A$ and $B$ given $C\cup D$. If there is no collider section on $\omega$ then there is a connecting walk between $A$ and $B$ given $C$, a contradiction. On $\omega$, all collider sections must have a node in $(C\cup D)$. If all collider sections have a node in $C$ then there is a connecting walk between $A$ and $B$ given $C$, again a contradiction. Hence consider first the collider section $\rho$ nearest  $A$ on $\omega$ that only has  nodes in $D$ on $\omega$; next, consider the closest node $i$ to $A$ on $\rho$ that  is in $D$. The subwalk between $A$ and $i$ then contradicts $A\dse B \cup D \cd C$.
%

4) \emph{Contraction:} If  $A\dse B\cd C$ and  $A\dse D\cd (B\cup C)$ then  $A\dse (B\cup D)\cd C$: Suppose, for contradiction, that there exists a connecting walk between $A$ and $B\cup D$ given $C$. Consider a shortest walk (i.e.\ a walk with fewest number of edges) of this type and call it $\omega$. The walk $\omega$ is either between $A$ and $B$ or between $A$ and $D$. The walk $\omega$ being between $A$ and $B$ contradicts $A\dse B\cd C$. Therefore, $\omega$ is between $A$ and $D$. In addition, since all collider sections on $\omega$ have a node in $C$ and  $A\dse D\cd (B\cup C)$, a non-collider section of $\omega$ must exist that has a node in $B\cup C$, and, therefore, in $B$. This contradicts the fact that $\omega$ is a shortest  connecting walk between $A$ and $B\cup D$ given $C$.

5) \emph{Intersection:} If $A\dse B\cd (C\cup D)$ and $A\dse D\cd (C\cup B)$ then $A\dse (B\cup D)\cd C$: Suppose, for contradiction, that there exists a connecting walk between $A$ and $B\cup D$ given $C$. Consider a shortest walk of this type and call it $\omega$. The walk $\omega$ is either between $A$ and $B$ or between $A$ and $D$. Because of symmetry between $B$ and $D$ in the formulation, it is enough to suppose that $\omega$ is between $A$ and $B$. Since all  collider sections on $\omega$ have a node in $C$ and  $A\dse B\cd (C\cup D)$, a non-collider section of $\omega$ must exist that has a node in $D\cup C$, and, therefore, in $D$. This contradicts the fact that $\omega$ is a shortest connecting walk between $A$ and $B\cup D$ given $C$.

6) \emph{Composition:} If $A\dse B\cd C$ and $A\dse D\cd C$ then $A\dse (B\cup D)\cd C$: Suppose, for contradiction, that there exist connecting walks between $A$ and $B\cup D$ given $C$. Consider a walk of this type and call it $\omega$. Walk $\omega$ is either between $A$ and $B$ or between $A$ and $D$. Because of symmetry between $B$ and $D$ in the formula it is enough to suppose that $\omega$ is between $A$ and $B$. But this contradicts $A\dse B\cd C$.
\end{proof}
This theorem implies that we can focus on establishing conditional independence for pairs of nodes, formulated in the corollary below.
\begin{coro}\label{cor:composition}
For a graph $G$ and disjoint subsets of nodes $A$, $B$, and $C$, it holds that $A\dse B\cd C$ if and only if $i\dse j\cd C$ for every pair of nodes $i\in A$ and $j\in B$.
\end{coro}
\begin{proof}
The result follows from the fact that $\dse$ satisfies  decomposition and  composition.
\end{proof}

\subsection{Relation to other separation criteria}
Four different types of independence models have previously been associated with chain graphs.  These are known as the LWF Markov property, defined by \cite{lau89} and later studied by e.g.\ \cite{fry90,stub98}; the AMP Markov property, defined and studied by \cite{and01}, and the multivariate regression (MR) Markov property, introduced by \cite{cox93} and studied e.g.\ by \cite{mar11}; in addition, \citet{drt09} briefly considered a type III chain graph Markov property which we are not further considering here.

Traditionally these have been formulated using undirected chain graphs but with different separation criteria.  In contrast, here we use a single notion of separation and the different independence models appear by varying the type of chain graph. In particular, the LWF Markov property corresponds to UCGs, the MR Markov property to BCGs, and the AMP Markov property to DCGs, as we shall see below.

Table~\ref{tab:1} gives an overview of different types of colliders used in the   various independence models associated with chain graphs.
\begin{table}[htb]
\caption{\small{Colliders for different chain graph independence models.}}\label{tab:1}
\vspace{5mm}
\centering
\begin{tabular}{|c|ccc|}
  \hline
  MR: & $\fra\circ\fla$ & $\fra\circ \arc$ & $\arc\circ\arc$ \\ \hline
  AMP: & $\fra\circ\fla$ & & $\fra\circ \ddash$  \\ \hline
  LWF: & & $\fra\circ\ful\cdots\ful\circ\fla$  &  \\
  \hline
\end{tabular}

\end{table}

For summary graphs and their subclasses, \cite{sadl14} showed that the unifying separation concept was that of \emph{$m$-separation}, defined as follows.
A path $\pi$  is \emph{$m$-connecting} given $C$ if all collider nodes on $\pi$ intersect $\An(C)$ and all non-collider nodes on $\pi$ are disjoint from $C$. Notice that $m$-separation considers nodes, but the fact that there is no arrowhead pointing to a node that is endpoint of a line in a summary graph implies that every collider section of any walk consists of a single node. For pairwise disjoint subsets $\langle A,B,C\rangle$, $A$ and $B$ are \emph{$m$-separated} by $C$ if there are no $m$-connecting paths between $A$ and $B$ given $C$, and we use the notation $A\mse B\cd C$ to indicate this.  The following lemma establishes that for summary graphs (and all subclasses of these), $m$-separation is equivalent to the separation we have defined here. The idea is similar to that employed in \cite{kos02}.
\begin{lemma}\label{lem:msepequiv}Suppose that $G$ is a summary graph. Then
$$A\dse B\cd C\iff A\mse B\cd C.$$
\end{lemma}
\begin{proof}
We need to show that for $i,j\not \in C$, there is a connecting walk between $i$ and $j$ if and only if there is an $m$-connecting path between $i$ and $j$ given $C$. If there is an $m$-connecting path $\pi$ between $i$ and $j$ then there exists a connecting walk between $i$ and $j$ by taking $\pi$ and add the possible directed path from a collider node $k$ on $\pi$ to $c\in C$ and its reverse from $c$ to $k$.

Thus suppose that there is a connecting walk $\omega$ between $i$ and $j$. Since there are no arrowheads pointing to nodes that are endpoints of lines, all collider sections on $\omega$ are single nodes; and hence we can talk of collider nodes instead of sections. Consider the walk between $i$ and $j$ obtained from $\omega$ by replacing any subwalk of type $\langle l, \rho', l\rangle$ (for a subwalk $\rho'$) by a single node subwalk $\langle l\rangle$. First of all, it is clear that the resulting walk is a path. Denote this path by $\pi$. We show that an $m$-connecting path can be constructed from $\pi$:

It is not possible that a node that occurs (at least once) as a collider on $\omega$ and occurs also as a member of a non-collider section on $\pi$: If $k$ is a collider node on $\omega$ then it is in $C$. This means that there is an arrowhead at $k$ on all tripaths with inner node $k$ on $\omega$. Hence, regardless of which two edges of $\omega$ with endpoint $k$ are on $\pi$, the corresponding tripath remains collider.

Therefore, all non-collider nodes on $\pi$ are outside $C$. If all collider nodes are in $C$ then we are done. Thus suppose that there is a collider node $k$ (on collider tripath $\langle k_0,k,k_1\rangle$) on $\pi$ that is not in $C$. This means that, on $\omega$, $k$ is always within a non-collider section. Consider an edge $kr_0$ on $\omega$ that is a part of the subwalk $\langle k_0,k,r_0\rangle$ of $\omega$, and notice that this edge is not on $\pi$. The edge $kr_0$ is not a line as otherwise there is an arrowhead pointing to
an endpoint of a line. As the edge $kr_0$ itself has no arrowhead at $k$ it must be
an arrow from $k$ to $r_0$. Following through $\omega$ from $r_0$, inductively, we have three cases: 1) There exists a directed cycle, which is impossible. 2) $k$ is an ancestor of a collider node $r$: We have that $r\in C$, and hence $k$ is an ancestor of $C$. 3) $k$ is an ancestor of $i$ or $j$: Without loss of generality, assume that $k\in\an(j)$. In this case, we modify $\pi$ by replacing the subwalk between $k$ and $j$ by a directed path from $k$ to $j$. Notice that no node on this path is in $C$. This completes the proof.
\end{proof}

For MAMPs, \cite{pen14} provides a generalization of the $p$-separation \cite{lev01} for AMP chain graphs. In the language and notations of this paper, it is defined a follows: A path $\pi$ is \emph{$z$-connecting}  given $C$ ($z$ is our notation) for MAMPs if every collider node on $\pi$ is in $\An(C)$ and every non-collider node $k$ is outside $C$ unless there is a subpath of $\pi$, $i\ddash k\ddash j$ such that  $\spo(k)\neq\varnothing$ or $\pa(k)\setminus C\neq\varnothing$. We say that $A$ and $B$ are $z$-separated given $C$, and write $A\dse_z B\cd C$, if there is no $z$-connecting path between $A$ and $B$ given $C$.
\begin{lemma}\label{lem:zsepequiv}Suppose that $G$ is a MAMP. Then
$$A\dse B\cd C\iff A\dse_z B\cd C.$$
\end{lemma}
\begin{proof}
We need to show that for $i,j\not \in C$, there is a connecting walk between $i$ and $j$ if and only if there is a $z$-connecting path between $i$ and $j$ given $C$. If there is a $z$-connecting path $\pi$ between $i$ and $j$ we may construct a connecting walk between $i$ and $j$ by modifying $\pi$ as follows: 1) for a collider node $k\in\an(C)$, add a directed path from $k$ to $c\in C$ and its reverse from $c$ to $k$; 2) for a non-collider node $k\in C$ within $i \ddash k \ddash j$ in $\pi$ (see the definition of $z$-separation), we distinguish two cases: if $\spo(k) \neq\varnothing$  then one has $i \ddash j$ by the definition of MAMP and one can shorten the tripath $i\ddash k \ddash j$ on $\pi$; if $\spo(k) = \varnothing$ but $l\in \pa(k)\setminus C$ exists then add the $kl$ edge and its reverse to $\pi$.
%

Thus suppose that there is a connecting walk $\omega$ between $i$ and $j$. Since there are no lines, all sections on $\omega$ are single nodes; and hence we can talk of collider and non-collider nodes instead of sections. Similar to Lemma \ref{lem:msepequiv}, consider the walk between $i$ and $j$ obtained from $\omega$, and whenever there is a node $l$ with repeated occurrence in $\omega$, replace the
cycle from $l$ to $l$ in $\omega$ by a single occurrence of $l$. The resulting walk is a path, denoted by $\pi$. We show that $z$-connecting path can be constructed from $\pi$:

The only case where a node $k$ is  a collider node on $\omega$ and it turns into a non-collider node on $\pi$ is when $k$ is the inner node of the tripath $h\ddash k\ddash l$ on $\pi$. Therefore, all non-collider nodes on $\pi$ are outside $C$ unless this mentioned case occurs. However, in this case either $\spo(k)\neq\varnothing$ or $\pa(k)\setminus C\neq\varnothing$, which ensures that the condition of the definition of a $z$-connecting path is still satisfied.

If all collider nodes are in $C$ then we are done.
Thus suppose that there is a collider node $k$ (on collider tripath $\langle k_0,k,k_1\rangle$) on $\pi$ that is not in $C$. This means that, on $\omega$, $k$ is always a non-collider node. There is an arrowhead at $k$ on at least one of the $k_0k$ or the $kk_1$ edges. Without loss of generality, assume that it is the $k_0k$ edge. Consider an edge $kr_0$ on $\omega$ that is a part of the subwalk $\langle k_0,k,r_0\rangle$ of $\omega$, and notice that this edge is not on $\pi$. As the edge $kr_0$ itself has no arrowhead at $k$ and is not a dotted line, it must be
an arrow from $k$ to $r_0$. Following through $\omega$ from $r_0$, inductively, we have three cases: 1) There exists a directed cycle, which is impossible. 2) $k$ is an ancestor of a collider node $r$: We have that $r\in C$, and hence $k$ is an ancestor of $C$. 3) $k$ is an ancestor of $i$ or $j$: Without loss of generality, assume that $k\in\an(j)$. In this case, we modify $\pi$ by replacing the subwalk between $k$ and $j$ by a directed path from $k$ to $j$. Notice that no node on this path is in $C$. This completes the proof.
\end{proof}

We are now ready to show that our concept of separation unifies the independence models discussed.
\begin{theorem} Independence models generated by separation in graphs with four types of edges are identical to the independence models associated with the subclasses in Fig.~\ref{fig:hierarchy}.
\end{theorem}
\begin{proof} It is shown in \cite{sadl14} that $m$-separation, as defined above, unifies independence models for SGs and subclasses thereof and by Lemma~\ref{lem:msepequiv} $m$-separation is equivalent to our separation. The separation criterion in \cite{sad16} for CMGs is identical to the separation given here when there are no dotted lines in the graph. Hence, the independence models generated by our separation criterion unifies independence models for all the subclasses of CMGs. Lemma \ref{lem:zsepequiv} shows that, dotted lines replacing lines in Pe\~na's separation criterion, it becomes identical to ours. For AADMGs, Criterion 2 defined as the global Markov property in \cite{pen16} is trivially a special case of the separation defined here. Therefore, our criterion unifies independence models in all subclasses of graphs.
\end{proof}
Notice that most of the associated classes of independence models presented in the diagram of Fig.~\ref{fig:hierarchy} are  distinct;  exceptions are AGs and SGs, which are alternative representations of the same class of independence models, and the same holds for DGs and UGs. In addition, it can be seen from Table \ref{tab:1} that, for every type of chain graph, one different type of symmetric edge is needed since each of them forms different colliders; hence, the unification for the general class of graphs with four types of edges is not achieved by graphs with three types of edges.
\subsection{Probabilistic independence models and the global Markov property}\label{sec:prob}
Consider a set $V$ and a collection of random variables
$(X_\alpha)_{\alpha\in V}$ with state spaces $\mathcal{X}_\alpha, \alpha\in V$ and joint distribution $P$. We let $X_A=(X_v)_{v\in A}$ etc.\ for each subset $A$ of $V$. For disjoint subsets $A$, $B$, and $C$ of $V$
we use the short notation $A\cip B\cd C$ to denote that $X_A$ is \emph{conditionally independent of $X_B$ given $X_C$} \citep{daw79,lau96}, i.e.\ that for any measurable $\Omega\subseteq \mathcal{X}_A$ and $P$-almost all $x_B$ and $x_C$,
$$P(X_A \in \Omega\cd X_B=x_B, X_C=x_C)=P(X_A \in \Omega\cd X_C=x_C).$$
We can now induce an independence model $\mathcal{J}(P)$ by letting
\begin{displaymath}
\langle A,B\cd C\rangle\in \mathcal{J}(P) \text{ if and only if } A\cip B\cd C \text{ w.r.t.\ $P$}.
\end{displaymath}

We note that for a probabilistic independence model $\mathcal{J}(P)$, the marginal independence model to a set $D=V\setminus M$ is the independence model generated by the marginal distribution. More formally, we define the \emph{marginal independence model} over a subset of the node set  $M$ as follows:
$$\alpha(\mathcal{J},M)=\{\langle A,B\cd C\rangle:\langle A,B\cd C\rangle\in\mathcal{J}\text{ and } (A\cup B\cup C)\cap M=\varnothing\},$$
which is defined over $V\setminus M$.

\begin{lemma}\label{lem:probmarg}
Let $\mathcal{J}(P)$ be a probabilistic independence model; its marginal independence model is the independence model generated by the marginal distribution, i.e.\ for $D=V\setminus M$ we have
$\alpha(\mathcal{J}(P),M)=\mathcal{J}(P_D).$
\end{lemma}
\begin{proof} This is immediate.\end{proof}

For a graph $G=(V,E)$, an independence model $\mathcal{J}$ defined over $V$ satisfies the \emph{global Markov property} w.r.t.\ a graph $G$, if for  disjoint subsets $A$, $B$, and $C$ of $V$ it holds that  $$A\dse B\cd C \implies \langle A,B\cd C\rangle\in \mathcal{J}.$$ 
If $\mathcal{J}(P)$ satisfies the global Markov property w.r.t.\ a graph $G$, we also say that \emph{$P$ is Markov w.r.t.\ $G$}.
We say that an independence model $\mathcal{J}$ is \emph{probabilistic} if there is a distribution $P$ such that  $\mathcal{J}= \mathcal{J}(P)$. We then also say that $P$ is \emph{faithful} to $\mathcal{J}$. If $P$ is faithful to $\mathcal{J}(G)$ for a graph $G$ then we also say that $P$ is \emph{faithful to  $G$}. Thus, if $P$ is faithful to $G$ it is also Markov w.r.t.\ $G$.

Probabilistic independence models are always semi-graphoids \citep{pea88}, whereas the converse is not necessarily true; see \citep{stu89}. If, for example, $P$ has strictly positive density, the induced independence model is always a graphoid; see e.g.\ Proposition 3.1 in \citep{lau96}. If the distribution $P$ is a regular multivariate Gaussian distribution, $\mathcal{J}(P)$ is a compositional graphoid; e.g.\ see \cite{stu05}.

Probabilistic independence models with positive densities are not in general compositional; this only holds for special types of multivariate distributions such as, e.g.\  Gaussian distributions and the symmetric binary distributions used in \cite{wer09}.
However, the following statement implies that it is not uncommon for a probabilistic independence model to satisfy composition:
\begin{prop}
If there is a graph $G$ to which $P$  is faithful,  then $\mathcal{J}(P)$ is a compositional graphoid.
\end{prop}
\begin{proof}
The result follows from Theorem \ref{thm:110} since then $\mathcal{J}(P)=\mathcal{J}(G)$.
\end{proof}

\section{Maximality for graphs}\label{sec:4}
A graph $G$ is called \emph{maximal} if adding an edge between any two non-adjacent nodes in $G$ changes the independence model $\mathcal{J}(G)$. Notice that in \cite{sadl14} the non-adjacency condition was incorrectly omitted.

\begin{figure}[htb]
\centering
\scalebox{0.25}{\includegraphics{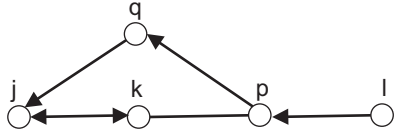}}
  \caption[]{\small{A non-maximal CMG.}}
     \label{fig:nonmaxex}
\end{figure}

Conditions 2 and 3, which MAMPs satisfy (provided in Section \ref{sec:2.1}) ensure that MAMPs are maximal; see \cite{pen14}.
However, graphs are not maximal in general.   For example, there exist non-maximal ancestral and summary graphs \cite{ric02,sadl14}; see also Fig.\ \ref{fig:nonmaxex} for an example of a graph that is neither a summary graph (hence it is not ancestral) nor maximal: this CMG induces no independence statements of the form $j\dse l \cd C$  for any choice of $C$: if we condition on $k$ or $p$ or both, the path $\langle j,k,p,l\rangle$ is connecting since $k\ful p$ is a collider section; conditioning on $q$ makes the walk $\langle j,k,p,q,p,l\rangle$ a connecting walk, and if we do not condition on anything, the walk  $\langle j,q,p,l\rangle$ is connecting.

\begin{figure}[htb]
\centering
\begin{tabular}{cc}
\scalebox{0.25}{\includegraphics{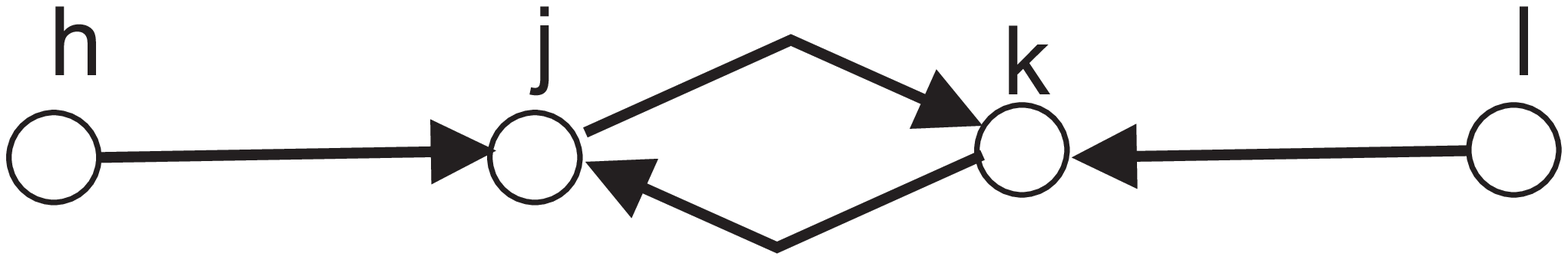}} &
\scalebox{0.25}{\includegraphics{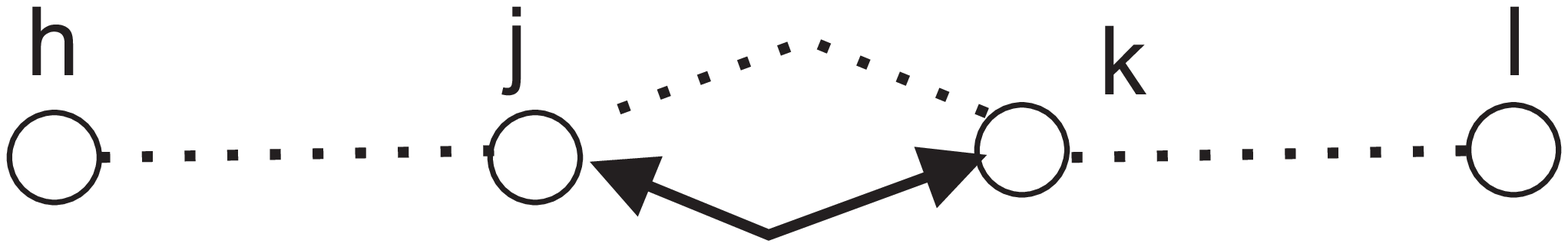}}\\
(a) & (b) 
\end{tabular}
  \caption[]{\small{Two non-maximal graphs that cannot be modified to be maximal by adding edges without changing the independence model.}}
     \label{fig:nonmaximizeable}
\end{figure}

The notion of maximality is important for pairwise Markov properties, to be discussed in the next section.
For a non-maximal ancestral or summary graph, one can obtain a maximal ancestral or summary graph with the same induced independence model by adding edges to the original graph \citep{ric02,sadl14}. As we shall show below, this is also true for general CMGs, but it is not generally the case for  graphs containing dotted lines or directed cycles. 
%
 Fig.~\ref{fig:nonmaximizeable} displays two small  \emph{non-maximizeable} graphs, where the graph in (a) contains a directed cycle.

For example, in the directed graph of Fig.~\ref{fig:nonmaximizeable}(a), in order to make the graph maximal, one must connect $h$ and $k$, and similarly $l$ and $j$. Now notice that in the original graph it holds that $h\dse l$ and $h\dse l\cd \{j,k\}$. However, after introducing new $hk$ and $lj$ edges, regardless of what type of edge we add, one of $h\dse l$ or $h\dse l\cd \{j,k\}$ does not hold. 

To characterise maximal CMGs we need the following notion:
A walk $\omega$  is a \emph{primitive inducing} walk between $i$ and $j$ ($i\neq j$) if and only if it is an $ij$ edge or $\omega= \langle i,q_1,q_2,\dots,q_p,j\rangle$ where for every $n$,
$1\leq n \leq p$, it holds that
\begin{description}
    \item[(i)] all inner sections of $\omega$ are colliders;
    \item[(ii)] endpoint sections of $\omega$ are single elements;
    \item[(iii)] $q_n\in \Ant(\{i,j\})$.
\end{description}
This definition is an extension of the notion of a primitive inducing path as defined for ancestral graphs in \citep{ric02}.  For example, in Fig.~\ref{fig:nonmaxex}, $\langle j, k, p,l\rangle$ is a primitive inducing walk. Next we need the following lemmas:
\begin{lemma}\label{lem:11300}
In a CMG, inner nodes of a walk $\omega$ between $i$ and $j$ that are on a non-collider section are either in $\ant(i)\cup\ant(j)$ or  anteriors of a collider section on $\omega$.
\end{lemma}
\begin{proof} Let $k=i_m$ be an inner node of $\omega$ and on a non-collider section on
a walk $\omega=\langle i=i_0,i_1,\dots,i_n=j\rangle$ in a CMG $G$. Then from at least one side
(say from $i_{m-1}$) there is no arrowhead on $\omega$ pointing to the section containing $k$. By moving towards $i$ on the path as long as $i_p$, $1\leq p\leq m-1$, is on a non-collider section on the walk, we obtain that $k\in\ant(i_{p-1})$.
This implies that if no $i_p$ is on a collider section then $k\in\ant(i)$.
\end{proof}
\begin{lemma}\label{lem:1130}
For nodes $i$ and $j$ in a CMG  that are not connected by any primitive inducing walks (and hence $i\not\sim j$), it holds that $i\dse j\cd \ant(\{i,j\})$.
\end{lemma}
\begin{proof}
Suppose that  there is a connecting walk $\varpi$ between $i$ and $j$ given $\ant(\{i,j\})$.

If $i$ or $j$ are on a non-collider inner section $\rho$ on $\varpi$ then $\rho$  is contained in $\{i, j\}$ since
otherwise any other node in $\rho$ would be in $\ant(\{i,j\})$, which is impossible. Then $\rho$ contains either only $i$ or only $j$ since $ij$
is not an edge in the graph. Thus, $\rho$  is either single $i$ or single $j$. In such a case remove the cycle between $i$ and $i$ (or between $j$ and $j$), which is a subwalk of $\varpi$. Repeat this process until there are no such non-collider sections. Denote the resulting walk by $\omega$. We shall show that $\omega$ is primitive inducing:

(i) If, for contradiction, there is a node $k$ on an inner non-collider section of $\omega$ then, by Lemma \ref{lem:11300}, $k$ is either in $\ant(i)\cup\ant(j)$  or it is an anterior of nodes of a collider section on $\omega$, but since $\omega$ is connecting given $\ant(\{i,j\})$, collider sections intersect $\ant(\{i,j\})$ and hence are in $\ant(\{i,j\})$ themselves. (Hence, $k\notin\{i,j\}$.) Now, $k\in\ant(\{i,j\})$ contradicts the fact that $\omega$ is connecting given $\ant(\{i,j\})$.  

(ii) Unless $\omega$ is a line, endpoint sections of $\omega$ are single elements since they are non-collider on $\omega$ and, if not single elements, their members, excluding $i$ or $j$, are in $\ant(\{i,j\})$, which is impossible.

(iii) This condition is clear since all inner nodes are in collider sections and consequently (except for possibly $i$ or $j$) in $\ant(i)\cup\ant(j)$.
\end{proof}
\begin{lemma}\label{lem:113n}
The only primitive inducing walk between $i$ and $j$ without
arrowheads at its endpoints (i.e.\ $i$ and $j$) is the line $ij$.
\end{lemma}
\begin{proof}
Consider such a walk $\omega$: Suppose, for contradiction, that there are other nodes other than $i,j$ on $\omega$, and assume that $q_1$ is the node adjacent to the endpoint $i$ on $\omega$ (i.e., there is $i\fra q_1$). Notice that $q_1\neq j$ since otherwise $i\fra j$, $j\fra q_p$, for some $p$, and $q_p\in\Ant(\{i,j\})$ lead to a contradiction.

Then the lack of semi-directed cycles implies that $q_1\in \ant (j)$ and hence there is another node $q_2$ on $\omega$. Similarly for $q_p$ adjacent to $j$ on $\omega$, $q_p\in \ant(i)$. But we may then construct a semi-directed cycle by taking the $iq_1$ edge, the anterior path from $q_1$ to $j$, the $jq_p$ edge, and the anterior path from $q_p$ back to $i$, a contradiction.
\end{proof}
Next we say that two
walks $\omega_1$ and $\omega_2$ (including edges) between $i$ and $j$ are  \emph{endpoint-identical} if there is  an arrowhead pointing to the endpoint section containing $i$ in $\omega_1$ if and only if there is an arrowhead pointing to the endpoint section containing $i$ in $\omega_2$ and similarly for $j$. For example, the paths $i\fra j$, $i\ful k\fra l\arc j$, and $i\fra k \arc l\ful j$ are all endpoint-identical as they have an arrowhead pointing to the section containing $j$ but  no arrowhead pointing to the section containing $i$ on the paths, but they are not endpoint-identical to $i\ful k\arc j$.
We then have the following:
\begin{lemma}\label{lem:113000}
Let $G$ be a CMG with the node set $V$. If there is a primitive inducing walk $\omega$ between $i$ and $j$ in $G$, and $C\subseteq V\setminus\{i,j\}$, then a
connecting walk between $i$ and $j$ given $C$ exists which is endpoint identical to $\omega$.
\end{lemma}
\begin{proof}
We denote the sections of the primitive inducing walk $\omega$ by $\langle i=\tau_0,\tau_1,\dots,\tau_{s-1},\tau_s=j\rangle$ and
note that if a section $\tau$ intersects $\ant(A)$ for any set $A$, it holds that $\tau\subseteq \ant(A)$. By Lemma \ref{lem:113n}, it is enough to consider two cases:

Case 1) There is an arrowhead at $j$ and no arrowhead at $i$ on $\omega$: First notice that the edge $iq_1$ is an arrow from $i$ to $q_1$. We construct an endpoint-identical connecting walk $\varpi$ given $C$ between $i$ and $j$. We start from $i$ and move towards $j$ on $\omega$ via $iq_1$ where $q_1\in\tau_1$. As long as along $\varpi$, a section $\tau_n$, $1\leq n\leq s-1$, intersects $\Ant(C)$, we do the following: If $\tau_n\cap C\neq \varnothing$ then we let $\varpi$ move to $\tau_{n+1}$. If $\tau_n\cap C =\varnothing$ but $\tau_n\subseteq \ant(C)$ then we let $\varpi$ move from $\tau_{n}$ to $C$ via an anterior path and back to $\tau_{n}$ by reversing this path, subsequently continuing to $\tau_{n+1}$ using the corresponding edge from $\omega$.

So suppose that possibly $\varpi$ reaches a section $\tau_m$ not intersecting $\Ant(C)$. Note that $\tau_m$ cannot only contain $i$ since otherwise it intersects $\Ant(C)$ (through $\tau_1$). If $\tau_m$ only contains $j$ then we already have a connecting walk. Hence, the only case that is left is when there is a $k\in\tau_m\setminus\{i,j\}$ such that $k\in\ant(\{i,j\})$. If $k\in\ant(i)$ then notice that $k$ is an anterior of $C$ through $i$ and $\tau_1$, which is impossible. Thus $k\in\ant(j)$ with no nodes on the anterior path in $C$. We can now complete $\varpi$ by letting it move to $j$ via this anterior path.

Notice that $\varpi$ is endpoint-identical to $\omega$ since both have an arrowhead at $j$ and no arrowhead at $i$.

Case 2) There is an arrowhead at $j$ and an arrowhead at $i$ on $\omega$: We follow the same method as in Case 1 to construct $\varpi$. The only difference is that $k\in\tau_m\setminus\{i,j\}$ can be in $\ant(i)$ without being an anterior of $C$. (In fact, $k$ and $i$ may be on the same section on $\omega$.) In this case we entirely replace the already constructed part of $\varpi$ by the reverse of the anterior path from $k$ to $i$ (which is from $i$ to $k$), and let $\varpi$ proceed to $\tau_{m+1}$.

Again it is clear that the constructed $\varpi$ and $\omega$  have an arrowhead at $j$. If $k$ and $i$ are not in the same section or are not connected by an undirected path then it is clear that there is an arrowhead at $i$, which is a single-node section on $\varpi$. If $k$ and $i$ are in the same section or are connected by an undirected path then there is an arrowhead at the endpoint section of $\varpi$ that contains $i$.

\end{proof}
Next, in Theorem~\ref{thm:112} we give a necessary and sufficient condition for a CMG  to be maximal. The analogous result for ancestral graphs was proved in Theorem 4.2 of \citep{ric02}.
\begin{theorem}\label{thm:112}
A CMG $G$ is maximal if and only if $G$ does not contain any primitive inducing walks between non-adjacent nodes.
\end{theorem}
\begin{proof}
($\Rightarrow$) Let $\omega =\langle i=i_0,i_1,\dots,i_n=j\rangle$ be a primitive inducing walk between non-adjacent nodes $i$ and $j$. By Lemma \ref{lem:113000}, there is therefore an endpoint-identical connecting walk $\omega'$ between $i$ and $j$ given any choice of $C$; thus, there is clearly no separation of form $i\dse j\cd C$. Let us add an endpoint-identical $ij$ edge to $G$. If a separation $A\dse B\cd C'$ is destroyed then the edge $ij$ is a part of the connecting walk $\omega''$ given $C'$ between $A$ and $B$. Now by replacing $ij$ by $\omega'$ on $\omega''$, we clearly obtain a walk in $G$ that is connecting given $C'$. This implies that adding $ij$ does not change $\mathcal{J}(G)$; hence, $G$ is not maximal.

($\Leftarrow$) By letting $C=\ant(\{i,j\})$ for every non-adjacent pair of nodes $i$ and $j$ and using Lemma \ref{lem:1130}, we conclude that for every missing edge there is an independence statement in $\mathcal{J}(G)$. This implies that $G$ is maximal.
\end{proof}
It now follows that for maximal graphs, every missing edge corresponds to a pairwise conditional independence statement in $\mathcal{J}(G)$:

\begin{coro}\label{prop:maxsep}
A CMG $G$ is maximal if and only if every missing edge in $G$ corresponds to a pairwise conditional independence statement in $\mathcal{J}(G)$.
\end{coro}
\begin{proof}
($\Leftarrow$) is clear. ($\Rightarrow$) follows from Theorem~\ref{thm:112} and Lemma \ref{lem:1130}.
\end{proof}
Also, we have the following corollary.
\begin{coro}
If $G$ is a non-maximal CMG then it can be made maximal by adding edges without changing its independence model.
\end{coro}
\begin{proof}
We begin with a non-maximal CMG $G$, and show that we can ``close" all the primitive inducing walks in order to obtain a maximal graph with the same induced independence model. For every primitive inducing walk $\omega$ between $i$ and $j$ where $i\nsim j$ in $G$, add an $ij$ edge that is endpoint-identical to $\omega$ if an edge of the same type does not already exist.

First we show that the resulting graph is a CMG: It is enough to show that an added edge does not generate a semi-directed cycle. By Lemma \ref{lem:113n}, the added edge is either an arrow or an arc. Since arcs are not on a semi-directed cycle, adding an arc would not generate a semi-directed cycle. Thus suppose that the added edge is an arrow from $i$ to $j$. Notice that the adjacent node $q_1$ to $i$ on the primitive inducing walk is in $\ant(j)$ and the $iq_1$ edge is an arrow from $i$ to $q_1$. Hence, if, for contradiction, the added $ij$ arrow generates a semi-directed cycle, a semi-directed cycle already existed in the original graph, where $ij$ is replaced by the anterior walk that consists of the $iq_1$ arrow and the anterior walk from $q_1$ to $j$. This is a contradiction.

Now, since the resulting graph does not contain any primitive inducing walks between non-adjacent nodes, it is maximal. In addition, by Lemma \ref{lem:113000}, there is a connecting walk between $i$ and $j$, which is endpoint-identical to the primitive inducing walk. One can replace the endpoint-identical $ij$ edge to this walk in any connecting walk in $G$ that contains $\omega$ as a subwalk.
\end{proof}

For example,  in Fig.~\ref{fig:nonmaxex}, $\langle j, k, p,l\rangle$ was a primitive inducing walk; hence this graph was not maximal. We may then add the edge $l\fra j$ and it becomes maximal.

\section{Pairwise Markov properties for chain mixed graphs}\label{sec:4n}
\subsection{A pairwise Markov property} It is possible to consider a general pairwise Markov property for specific subclasses of graphs with four types of edges (that actually have the four types) by including the results of  \cite{pen14,pen16}, which define  pairwise Markov properties for  marginal AMP chain graphs and alternative ADMGs and show the equivalence of pairwise and global Markov properties for such graphs. However, such a unification would be technically complex. Hence, we henceforth focus on CMGs; thus the considerations here concerning pairwise Markov properties do not cover AMP chain graphs.
%

A pairwise Markov property provides independence statements for non-adjacent pairs of nodes in the graph. For maximal graphs any non-adjacent nodes $i$ and $j$ are independent given some set $C$, but a pairwise Markov property yields a specific choice of $S=S(i,j)$ for every non-adjacent pair $i,j$.  The choice we provide here for any CMG immediately extends the choice in \cite{sadl14}. We show that for a maximal CMG, this pairwise Markov property is equivalent to the global Markov property for compositional graphoid independence models; in other words,  the pairwise statements combined with the compositional graphoid axioms generate the full independence model.  The maximality is critical for the pairwise statements to hold, as discussed above.

%
%

An independence model $\mathcal{J}$ defined over $V$ satisfies the \emph{pairwise Markov property} (P) w.r.t.\ a CMG $G$ if for every pair of nodes $i$ and $j$ with $i\not\sim j$ it holds that 
$$\mbox{(P)}:\nn\nn\nn  \langle i,j\cd \ant(\{i,j\})\rangle\in \mathcal{J}.$$

The pairwise Markov property simplifies for specific subclasses of graphs.
For connected UGs we have $\ant(\{i,j\})=V\setminus\{i,j\}$ and hence the standard pairwise Markov property appears; and for BGs we have $\ant(\{i,j\})=\varnothing$, so the property is identical to pairwise independence of non-adjacent nodes. For SGs and AGs (which include DAGs), a semi-direction preserving path is of the form $\circ\ful\dots\circ\ful\circ\fra\dots\circ\fra$, hence the anterior path (and consequently (P)) specializes to those in \cite{sadl14} and \cite{ric02} respectively.

Strictly speaking,  the unification only contains ``connected" UGs. It is not possible to extend the unification to all UGs and at the same time keep the pairwise Markov properties defined in the literature for other classes under any unified pairwise Markov property: In principle, it is fine to add nodes that are not in the connected component(s) of $i$ and $j$ to the conditioning set in any pairwise Markov property. However, although the well-known pairwise Markov property for UGs contains all such nodes, the known pairwise Markov properties for other classes do not.



\subsection{Equivalence of pairwise and global Markov properties}
Before establishing the main result of this section, we need several lemmas. We shall need to consider marginalization of independence models and use that it preserves the compositional graphoid property, shown in Lemma 8 of \cite{sadl14}:
\begin{lemma}\label{lem:j1}
Let $\mathcal{J}$ be a compositional graphoid over a set $V$ and $M$ a subset of $V$. It then holds that the marginal independence model $\alpha(\mathcal{J},M)$ is also a compositional graphoid.
\end{lemma}
Moreover,  we have
\begin{lemma}\label{lem:anterior} Let $\mathcal{J}=\mathcal{J}(G)$ be the independence model induced by a CMG $G$ and $M\subseteq V$. If $D=V\setminus M$ is an anterior set, the marginal model is determined by the induced subgraph $G[D]$:
$$\alpha(\mathcal{J}(G),M)= \mathcal{J}(G[D]).$$
\end{lemma}
\begin{proof}  We need to show that for $\{i,j\}\cup C\subseteq D$ we have that $i\dse j\cd C$ if and only if this is true in the induced subgraph $G[D]$. Clearly, if a connecting walk between $i$ and $j$ runs entirely within $D$ it also connects in $G$. Assume for contradiction that there is a connecting walk  which has a node  $k$ outside $D$ and consider an excursion on the walk that leaves $D$ at $i^*$, reaches $k$, and reenters into $D$ at $j^*$.  Since the walk is connecting, there are no collider sections on this excursion and thus it follows from Lemma~\ref{lem:11300} that $k$ is either anterior to $i^*$ or to $j^*$,  which contradicts the fact that $D$ is an anterior set.
\end{proof}

The following important lemma and its corollary imply that for any non-adjacent pair $i\not\sim j$ in a maximal CMG we can always find an $(i,j)$-separator $C$ with  $C\subseteq \ant(\{i,j\})$.
\begin{lemma}\label{lem:j3}
For a pair of distinct nodes $i$ and $j$ and a subset of the node set $C$ in a maximal CMG, if $i\dse  j\cd C$ for $C\setminus\ant(\{i,j\})\neq\varnothing$, then there is a node $l\not\in\ant(\{i,j\})$ in $C$ such that $i\dse  j\cd C\setminus\{l\}$.
\end{lemma}
\begin{proof}
Let $l'\in C\setminus\ant(\{i,j\})$ be arbitrary. If there is an $l''\in C\setminus\ant(\{i,j\})$ so that $l'\in\ant(l'')$ but $l''\notin\ant(l')$, then replace $l'$ by $l''$, and repeat this process until it terminates, which is ensured by the transitivity of semi-directed walks and the lack of semi-directed cycles in the CMG. Call the final node $l$. Thus, if $l\in\ant(\tilde{l})$ for $\tilde{l}\in C\setminus\ant(\{i,j\})$ then we also have that $\tilde{l}\in\ant(l)$. The lack of semi-directed cycles implies that this is equivalent to $l$ and $\tilde{l}$ being connected by lines.

We now claim that $i\dse j\cd C\setminus \{l\}$.
Suppose, for contradiction, that there is a connecting walk $\omega$ between $i$ and $j$ given $C\setminus\{l\}$. If $l$ is not on $\omega$ then $\omega$ is also connecting given $C$.
In addition, we have that $l$ is on a non-collider section $\rho$ on $\omega$. There is no arrowhead at $\rho$ from at least one side of the section, say from the $i$ side. We move towards $i$ on $\omega$ and denote the corresponding subwalk of $\omega$ by $\omega'=\langle l=l_0,l_1,\dots,l_m=i\rangle$. As long as $l_p$, $1\leq p\leq m-1$, is on a non-collider section on $\omega'$, we obtain that there is a semi-directed walk from $l$ to $l_p$. This implies that if no $l_p$ is on a collider section then there is an anterior walk from $l$ to $i$, which is impossible.

Therefore, by moving towards $i$ from $l$, we first reach an $\tilde l$ on $\omega'$ that lies on a collider section and is in $C\setminus\{l\}$.  Transitivity of anterior walks and the fact that there is no anterior walk from $l$ to $i$ or $j$ now imply that there is no anterior walk from $\tilde l$ to $i$ or $j$. The construction of $l$  implies that $l$ and $\tilde l$ are on the same section, and hence $l$ is not on a non-collider section on $\omega$, a contradiction.
Hence we conclude that $i\dse j\cd C\setminus \{l\}$.
\end{proof}
\begin{coro}\label{cor:antpair}
For a pair of nodes $i$ and $j$ and a subset $C$ of the node set in a maximal CMG, if $i\dse  j\cd C$, then $i\dse  j\cd C\cap \ant(\{i,j\})$.
\end{coro}
\begin{proof}Lemma~\ref{lem:j3} implies that we can repeatedly remove single nodes in $C\setminus \ant(\{i,j\})$ and preserve separation to obtain that  $i\dse j\cd C\cap \ant(\{i,j\})$. This concludes the proof.
\end{proof}
A direct implication of Lemma \ref{lem:1130} and Theorem \ref{thm:112} establishes  that the induced independence model $\mathcal{J}(G)$ for a maximal CMG $G$ satisfies the pairwise Markov property (P):
\begin{prop}\label{prop:g_implies_p}
If $i\not\sim j$ are non-adjacent nodes in a maximal CMG $G$, it holds that $i\dse j\cd \ant(\{i,j\})$.
\end{prop}
Finally we are ready to show the main result of this section.
\begin{theorem}\label{thm:114}
Let $G$ be a maximal CMG. If an independence model $\mathcal{J}$  over the node set of $G$ is a compositional graphoid, then $\mathcal{J}$ satisfies the pairwise Markov property (P) w.r.t.\ $G$ if and only if it satisfies the global Markov property w.r.t.\ $G$.
\end{theorem}
\begin{proof}

That the global Markov property implies the pairwise property (P) follows directly from Proposition~\ref{prop:g_implies_p}.

Now suppose that $\mathcal{J}$ satisfies the pairwise Markov property (P) and compositional graphoid axioms.  For subsets $A$, $B$,
and $C$ of the node set of $G$, we must show that $A\dse  B\cd C$ implies $\langle A,B\cd C\rangle\in \mathcal{J}$.
By Corollary~\ref{cor:composition}, it is
sufficient to show this when $A$ and $B$ are singletons, i.e.\ that $i\dse  j\cd C$ implies $\langle i,j\cd C\rangle\in \mathcal{J}$.

We establish the result in two main parts. In part I we consider the case with $C\subseteq \ant(\{i,j\})$ and  in part II we consider the general
case.

\subparagraph{Part I}  Suppose that $C\subseteq \ant(\{i,j\})$. We use induction on the number of nodes of the graph. The induction base
for a graph with two nodes is trivial. Thus suppose that the conclusion holds for all graphs with fewer than $n$ nodes and assume that $G$
has $n$ nodes.

Suppose there is an anterior set $D$ such that  $M=V\setminus D\neq \varnothing$ and $\{i\}\cup\{j\}\cup C\subseteq D$. The marginal independence model $\alpha(\mathcal{J},M)$ clearly also satisfies the pairwise Markov property w.r.t.\ $G[D]$ and hence the inductive assumption together
with Lemmas \ref{lem:j1} and \ref{lem:anterior} yields $\langle\{i\},\{j\}\cd C\rangle\in \mathcal{J}$. 
%
%
%

So suppose that this is not the case and hence $V=Ant(\{i,j\})$.
We establish the conclusion by reverse induction on $|C|$: For
the base we have $C=V\setminus\{i,j\}=\ant(\{i,j\})$ and the result follows directly from the pairwise Markov property.

For the inductive step, consider a node $h\not\in C$. We want to show that there are not simultaneously connecting walks between $h$ and $i$ and $h$ and $j$:  Suppose, for contradiction, there are connecting walks $\omega_1=\langle i,i_1,\dots,i_n,h\rangle$ and
$\omega_2=\langle h,j_m,j_{m-1},\dots,j_0=j\rangle$ given $C$. If, on the walk $\langle\omega_1,\omega_2\rangle$, the node $h$ is on a non-collider section then  so is $h$ on both $\omega_1$ and $\omega_2$, and hence $i$ and $j$ are connected given $C$, a contradiction. Thus we need only consider the case where $h$ is on a collider section on $\langle\omega_1,\omega_2\rangle$. However, we know that $h\in\ant(i)$ or $h\in\ant(j)$. Because of symmetry between $i$ and $j$ suppose that $h\in\ant(j)$, and denote the anterior path from $h$ to $j$ by $\omega_3$. Notice that the section containing $h$ on $\omega_1$ is non-collider and hence all members are outside $C$. Now, if no node on $\omega_3$ is in $C$ then $\langle\omega_1,\omega_3\rangle$ is a connecting walk between $i$ and $j$ a contradiction; and if there is a node $k$ on $\omega_3$ is in $C$ then $\langle\omega_1,\omega_4,\omega_4^r,\omega_2\rangle$ is a connecting walk between $i$ and $j$, where $\omega_4$ is the subwalk of $\omega_3$ between $h$ and $k$ and $\omega_4^r$ is $\omega_4$ in reverse direction, a contradiction again.
We conclude that, given $C$, $h$ is not connected to both $i$ and $j$.

By symmetry suppose that $i\dse  h\cd C$.
We also have that $i\dse  j\cd C$. Since $\mathcal{J}(G)$ is a compositional graphoid (Theorem~\ref{thm:110}) the composition property gives
that $i\dse  \{j,h\}\cd C$. By weak union for $\dse $ we obtain $i\dse  j\cd \{h\}\cup C$ and $i\dse  h\cd \{j\}\cup C$. By the induction hypothesis
we obtain $\langle i,j\cd \{h\}\cup C\rangle\in \mathcal{J}$ and $\langle i,h\cd \{j\}\cup C\rangle\in \mathcal{J}$. By intersection we get
$\langle i,\{j,h\}\cd C\rangle\in \mathcal{J}$. By decomposition we finally obtain $\langle i,j\cd C\rangle\in \mathcal{J}$.

\subparagraph{Part II} We now prove the result in the general case by induction on $|C|$. The base, i.e.\ the case that $|C|=0$, follows from part I. To prove
the inductive step we can assume that $C\nsubseteq \ant(\{i,j\})$, since otherwise part I implies the result.

By Lemma \ref{lem:j3}, since $C\nsubseteq \ant(\{i,j\})$, there is a node $l\in C$ such that $i\dse  j\cd C\setminus\{l\}$. We now have that either $i\dse  l\cd C\setminus\{l\}$ or $j\dse  l\cd C\setminus\{l\}$
since otherwise there is a connecting walk between $i$ and $j$ given $C\setminus\{l\}$ in the case that $l$ is on a
non-collider section or given $C$ in the case that $l$ is on a  collider section. Because of symmetry, suppose that $i\dse  l\cd C\setminus\{l\}$.
By the induction hypothesis we have $\langle i,j\cd C\setminus\{l\}\rangle\in \mathcal{J}$ and $\langle i,l\cd C\setminus\{l\}\rangle\in \mathcal{J}$.
By the composition property we get $\langle i,\{j, l\}\cd C\setminus\{l\}\rangle\in \mathcal{J}$. The weak union property implies
$\langle i,j\cd C\rangle\in \mathcal{J}$.
\end{proof}
If we specialize Theorem~\ref{thm:114} to the most common case of probabilistic independence models, we get:
\begin{coro}
Let $G$ be a maximal CMG. A probabilistic independence model that satisfies the intersection and composition axioms satisfies the pairwise Markov property (P) w.r.t.\ $G$ if and only if it satisfies the global Markov property w.r.t.\ $G$.
\end{coro}

The theorem states that the intersection and composition properties are sufficient for equivalence of pairwise and global Markov properties. Notice that they are also necessary since for example for the simple subclass of DAGs they are also necessary; see Section 6.3 of \cite{sadl14}.

\subsection{Alternative pairwise Markov properties} There are typically many other valid choices of the separating sets  $C(i,j)$ defining the pairwise Markov properties, see for example \cite{pen14}. In general a pairwise Markov property (P*) has the form
\begin{equation*}
  \mbox{(P*)}: \nn\nn i\nsim j\nn\Rightarrow\nn
                     \langle i,j\cd C(i,j)\rangle\in \mathcal{J},
\end{equation*}
where $C(i,j)$ is an $(i,j)$-separator in $G$ for every $(i,j)$. The question then is what are the possible choices of such separator systems which would ensure that these separations form a `basis' for the independence model $\mathcal{J}$ in the sense that all conditional independences in $\mathcal{J}$ can be derived from (P*) using the compositional graphoid axioms. The example below shows that not all choices of separator systems are possible.
\begin{example}\rm \label{ex:counter}Consider the independence model $\mathcal{J}$ over $V=\{ 1,2,3,4,5\}$ containing the statements
\[\langle1,3\cd 2\rangle, \langle1,4\cd 3\rangle,\langle1,5\cd 4\rangle,\langle2,4\cd 1,3,5\rangle,\langle2,5\cd 3\rangle,\langle3,5\cd 1,2,4\rangle\]
as well as their symmetric counterparts and all independence statements of the form $\langle A, \emptyset\cd B\rangle$ or $\langle\emptyset, A\cd B\rangle.$
This independence model is easily seen to satisfy the compositional graphoid axioms. In addition, if we let $G$ be the graph $1\ful2\ful3\ful4\ful5$, each of the conditioning sets for statements of the form $\langle i,j\cd C(i,j)\rangle$ in $\mathcal{J}$ are indeed $(i,j)$-separators in $G$. Thus $\mathcal{J}$ satisfies (P*) w.r.t.\ the graph $G$, but clearly it does not satisfy the global Markov property w.r.t.\ $G$.
\qed
\end{example}
Note that the independence model in Example~\ref{ex:counter} may not be probabilistically representable. It is unclear to us whether the pairwise statements in (P*) for any system $C(i,j)$ of $(i,j)$-separators in an undirected graph, say, is sufficient to generate all independence statements of the form $\langle A,B\cd C\rangle$ for a \emph{probabilistic} compositional graphoid. The standard probabilistic counterexamples \cite{matus:92} involving the pairwise Markov properties are not compositional graphoids, hence are not counterexamples in this context.

For the subclass of regression graphs, four different pairwise Markov properties were defined in \cite{sadw15}, which are all equivalent to the global Markov property and to each other under compositional graphoid axioms.
\section{Summary and conclusion}\label{sec:5n}
In this paper, we used a similar approach to that of \cite{sadl14} to unify Markov properties for most classes of graphs in the literature of graphical models. 
The general idea is that for any of the three standard interpretations of the chain graph Markov property, (LWF, AMP, and multivariate regression), we use one type of edge in the unifying class of graphs and then use  a single separation criterion which is a natural generalization of $c$-separation as defined in \cite{stub98}.

Unifying an equivalent pairwise Markov property seems very technical when including the AMP chain graphs, hence we restricted ourselves to prove the equivalence of pairwise and global Markov properties for the class of maximal CMGs, which includes connected chain graphs with the LWF interpretation  as well as maximal summary graphs (and consequently maximal ancestral graphs). In order to include the class of AMP chain graphs or its generalizations for the unification of the pairwise Markov property, excluding certain ``directed cycles'' from the class of graphs with four types of edges (similar to the exclusion of semi-directed cycles in mixed graphs) is necessary.

It was seen in this paper that, under compositional graphoid axioms, the system of pairwise independence statements constituting the pairwise Markov property, can act as a generating class for all independence statements given by the global Markov property. Typically there are many other systems of pairwise statements  that may act as a generating class for the global Markov property. The point given here is that there is a unified choice of these statements for the case of CMGs.

The two important independence models are induced by graphs and probability distributions. Establishing the pairwise Markov property for the independence model induced by graphs suffices for establishing the global Markov property as it is always a compositional graphoid. This is not always the case for the independence model induced by any probability distribution as the intersection and composition properties may not hold for such distributions.

\section*{Acknowledgements}
The authors are very grateful to the anonymous reviewers for careful, detailed, and helpful comments.
\bibliographystyle{imsart-nameyear}
\bibliography{bib}

\end{document}

%% file: deff3.tex
\usepackage{graphicx}

\usepackage{here}

\usepackage{latexsym}

\usepackage{amsbsy}

\usepackage{amsmath,amsthm,amssymb, amsfonts}

\usepackage{dashrule}

\newcommand{\pa}{\mathrm{pa}}
\newcommand{\nei}{\mathrm{ne}}
\newcommand{\spo}{\mathrm{sp}}
\newcommand{\an}{\mathrm{an}}
\newcommand{\An}{\mathrm{An}}
\newcommand{\partn}{\mathrm{pt}}

\newcommand{\ant}{\mathrm{ant}}
\newcommand{\Ant}{\mathrm{Ant}}

\newcommand{\dse}{\,\mbox{$\perp$}\,}

\newcommand{\mse}{\,\mbox{$\perp_{m}$}\,}
\newcommand{\cip}{\mbox{\,$\perp\!\!\!\perp$\,}}

\newcommand{\notdse}{\nolinebreak{\not\hspace{-1.5mm}\dse}}

\newcommand{\cd}{\,|\,}

\newcommand{\nn}[0]{\hspace*{.7em}}

\newcommand{\ful}{\mbox{$\, \frac{ \nn \nn \;}{ \nn \nn
}\,$}}

\newcommand{\fla}{\mbox{$\hspace{.05em} \prec
\!\!\!\!\!\frac{\nn \nn}{\nn}$}}

\newcommand{\fra}{\mbox{$\hspace{.05em} \frac{\nn
\nn}{\nn
}\!\!\!\!\! \succ \! \hspace{.25ex}$}}

\newcommand{\arc}{\mbox{$\hspace{.06em} \prec
\!\!\!\!\!\frac{\nn \nn}{\nn}
\!\!\!\!\!
\succ\! \hspace{.25ex}$}}

\newcommand{\ddash}{
{\,\cdot\!\cdot\!\cdot\!\cdot\!\cdot\,}}

\newtheorem{prop}{Proposition}
\newtheorem{coro}{Corollary}

\newtheorem{lemma}{Lemma}

\newtheorem{theorem}{Theorem}
\newtheorem{example}{Example}

%% file: CMGs-arxiv.bbl
\begin{thebibliography}{43}

\bibitem[\protect\citeauthoryear{Andersson, Madigan and Perlman}{2001}]{and01}
\begin{barticle}[author]
\bauthor{\bsnm{Andersson},~\bfnm{Steen~A.}\binits{S.~A.}},
  \bauthor{\bsnm{Madigan},~\bfnm{David}\binits{D.}} \AND
  \bauthor{\bsnm{Perlman},~\bfnm{Michael~D.}\binits{M.~D.}}
(\byear{2001}).
\btitle{Alternative {M}arkov properties for chain graphs}.
\bjournal{Scand.\ J.\ Stat.}
\bvolume{28}
\bpages{33--85}.
\end{barticle}
\endbibitem

\bibitem[\protect\citeauthoryear{Cox and Wermuth}{1993}]{cox93}
\begin{barticle}[author]
\bauthor{\bsnm{Cox},~\bfnm{D.~R.}\binits{D.~R.}} \AND
  \bauthor{\bsnm{Wermuth},~\bfnm{N.}\binits{N.}}
(\byear{1993}).
\btitle{Linear dependencies represented by chain graphs (with discussion)}.
\bjournal{Statist.\ Sci..}
\bvolume{8}
\bpages{204--218; 247--277}.
\end{barticle}
\endbibitem

\bibitem[\protect\citeauthoryear{Darroch, Lauritzen and Speed}{1980}]{dar80}
\begin{barticle}[author]
\bauthor{\bsnm{Darroch},~\bfnm{J.~N.}\binits{J.~N.}},
  \bauthor{\bsnm{Lauritzen},~\bfnm{S.~L.}\binits{S.~L.}} \AND
  \bauthor{\bsnm{Speed},~\bfnm{T.~P.}\binits{T.~P.}}
(\byear{1980}).
\btitle{Markov fields and log-linear interaction models for contingency
  tables}.
\bjournal{Ann.\ Statist.}
\bvolume{8}
\bpages{522--539}.
\end{barticle}
\endbibitem

\bibitem[\protect\citeauthoryear{Dawid}{1979}]{daw79}
\begin{barticle}[author]
\bauthor{\bsnm{Dawid},~\bfnm{A.~P.}\binits{A.~P.}}
(\byear{1979}).
\btitle{Conditional independence in statistical theory (with discussion)}.
\bjournal{J.\ Roy.\ Statist.\ Soc.\ Ser.\ B}
\bvolume{41}
\bpages{1--31}.
\end{barticle}
\endbibitem

\bibitem[\protect\citeauthoryear{Didelez}{2008}]{didelez:08}
\begin{barticle}[author]
\bauthor{\bsnm{Didelez},~\bfnm{Vanessa}\binits{V.}}
(\byear{2008}).
\btitle{Graphical models for marked point processes based on local
  independence}.
\bjournal{J.\ Roy.\ Statist.\ Soc.\ Ser.\ B}
\bvolume{70}
\bpages{245--264}.
\bdoi{10.1111/j.1467-9868.2007.00634.x}
\end{barticle}
\endbibitem

\bibitem[\protect\citeauthoryear{Drton}{2009}]{drt09}
\begin{barticle}[author]
\bauthor{\bsnm{Drton},~\bfnm{M.}\binits{M.}}
(\byear{2009}).
\btitle{Discrete chain graph models}.
\bjournal{Bernoulli}
\bvolume{15}
\bpages{736--753}.
\end{barticle}
\endbibitem

\bibitem[\protect\citeauthoryear{Drton and Richardson}{2008}]{drt08}
\begin{barticle}[author]
\bauthor{\bsnm{Drton},~\bfnm{Mathias}\binits{M.}} \AND
  \bauthor{\bsnm{Richardson},~\bfnm{Thomas~S.}\binits{T.~S.}}
(\byear{2008}).
\btitle{Binary models for marginal independence}.
\bjournal{J.\ Roy.\ Statist.\ Soc.\ Ser.\ B}
\bvolume{41}
\bpages{287--309}.
\end{barticle}
\endbibitem

\bibitem[\protect\citeauthoryear{Eichler}{2007}]{eichler:07}
\begin{barticle}[author]
\bauthor{\bsnm{Eichler},~\bfnm{Michael}\binits{M.}}
(\byear{2007}).
\btitle{Granger causality and path diagrams for multivariate time series}.
\bjournal{J.\ Econometrics}
\bvolume{137}
\bpages{334--353}.
\bdoi{http://dx.doi.org/10.1016/j.jeconom.2005.06.032}
\end{barticle}
\endbibitem

\bibitem[\protect\citeauthoryear{Frydenberg}{1990}]{fry90}
\begin{barticle}[author]
\bauthor{\bsnm{Frydenberg},~\bfnm{M.}\binits{M.}}
(\byear{1990}).
\btitle{The chain graph {M}arkov property}.
\bjournal{Scand. J. Stat.}
\bvolume{17}
\bpages{333--353}.
\end{barticle}
\endbibitem

\bibitem[\protect\citeauthoryear{Geiger, Verma and Pearl}{1990}]{gei90}
\begin{barticle}[author]
\bauthor{\bsnm{Geiger},~\bfnm{D.}\binits{D.}},
  \bauthor{\bsnm{Verma},~\bfnm{T.~S.}\binits{T.~S.}} \AND
  \bauthor{\bsnm{Pearl},~\bfnm{J.}\binits{J.}}
(\byear{1990}).
\btitle{Identifying independence in {B}ayesian networks}.
\bjournal{Networks}
\bvolume{20}
\bpages{507--534}.
\end{barticle}
\endbibitem

\bibitem[\protect\citeauthoryear{Gibbs}{1902}]{gibbs:02}
\begin{bbook}[author]
\bauthor{\bsnm{Gibbs},~\bfnm{W.}\binits{W.}}
(\byear{1902}).
\btitle{Elementary Principles of Statistical Mechanics}.
\bpublisher{Yale University Press}, \baddress{NewHaven, Connecticut}.
\end{bbook}
\endbibitem

\bibitem[\protect\citeauthoryear{Kauermann}{1996}]{kau96}
\begin{barticle}[author]
\bauthor{\bsnm{Kauermann},~\bfnm{G.}\binits{G.}}
(\byear{1996}).
\btitle{On a dualization of graphical {G}aussian models}.
\bjournal{Scand.\ J.\ Stat.}
\bvolume{23}
\bpages{105--116}.
\end{barticle}
\endbibitem

\bibitem[\protect\citeauthoryear{Kiiveri, Speed and Carlin}{1984}]{kii84}
\begin{barticle}[author]
\bauthor{\bsnm{Kiiveri},~\bfnm{H.}\binits{H.}},
  \bauthor{\bsnm{Speed},~\bfnm{T.~P.}\binits{T.~P.}} \AND
  \bauthor{\bsnm{Carlin},~\bfnm{J.~B.}\binits{J.~B.}}
(\byear{1984}).
\btitle{Recursive causal models}.
\bjournal{J.\ Aust.\ Math.\ Soc.\ Ser.\ A}
\bvolume{36}
\bpages{30--52}.
\end{barticle}
\endbibitem

\bibitem[\protect\citeauthoryear{Koster}{1996}]{koster:96}
\begin{barticle}[author]
\bauthor{\bsnm{Koster},~\bfnm{J.~T.~A.}\binits{J.~T.~A.}}
(\byear{1996}).
\btitle{Markov properties of nonrecursive causal models}.
\bjournal{Ann.\ Statist.}
\bvolume{24}
\bpages{2148--2177}.
\end{barticle}
\endbibitem

\bibitem[\protect\citeauthoryear{Koster}{2002}]{kos02}
\begin{barticle}[author]
\bauthor{\bsnm{Koster},~\bfnm{J.~T.~A.}\binits{J.~T.~A.}}
(\byear{2002}).
\btitle{Marginalizing and conditioning in graphical models}.
\bjournal{Bernoulli}
\bvolume{8}
\bpages{817--840}.
\end{barticle}
\endbibitem

\bibitem[\protect\citeauthoryear{Lauritzen}{1996}]{lau96}
\begin{bbook}[author]
\bauthor{\bsnm{Lauritzen},~\bfnm{S.~L.}\binits{S.~L.}}
(\byear{1996}).
\btitle{Graphical Models}.
\bpublisher{Clarendon Press}, \baddress{Oxford, United Kingdom}.
\end{bbook}
\endbibitem

\bibitem[\protect\citeauthoryear{Lauritzen and
  Spiegelhalter}{1988}]{lauritzen:spiegelhalter:88}
\begin{barticle}[author]
\bauthor{\bsnm{Lauritzen},~\bfnm{S.~L.}\binits{S.~L.}} \AND
  \bauthor{\bsnm{Spiegelhalter},~\bfnm{D.~J.}\binits{D.~J.}}
(\byear{1988}).
\btitle{Local computations with probabilities on graphical structures and their
  application to expert systems (with discussion)}.
\bjournal{J.\ Roy.\ Statist.\ Soc.\ Ser.\ B}
\bvolume{50}
\bpages{157--224}.
\end{barticle}
\endbibitem

\bibitem[\protect\citeauthoryear{Lauritzen and Wermuth}{1989}]{lau89}
\begin{barticle}[author]
\bauthor{\bsnm{Lauritzen},~\bfnm{S.~L.}\binits{S.~L.}} \AND
  \bauthor{\bsnm{Wermuth},~\bfnm{N.}\binits{N.}}
(\byear{1989}).
\btitle{Graphical models for association between variables, some of which are
  qualitative and some quantitative}.
\bjournal{Ann. Statist.}
\bvolume{17}
\bpages{31--57}.
\end{barticle}
\endbibitem

\bibitem[\protect\citeauthoryear{Levitz, Perlman and Madigan}{2001}]{lev01}
\begin{barticle}[author]
\bauthor{\bsnm{Levitz},~\bfnm{Michael}\binits{M.}},
  \bauthor{\bsnm{Perlman},~\bfnm{Michael~D.}\binits{M.~D.}} \AND
  \bauthor{\bsnm{Madigan},~\bfnm{David}\binits{D.}}
(\byear{2001}).
\btitle{Separation and completeness properties for amp chain graph Markov
  models}.
\bjournal{Ann. Statist.}
\bvolume{29}
\bpages{1751--1784}.
\end{barticle}
\endbibitem

\bibitem[\protect\citeauthoryear{Marchetti and Lupparelli}{2011}]{mar11}
\begin{barticle}[author]
\bauthor{\bsnm{Marchetti},~\bfnm{Giovanni~M.}\binits{G.~M.}} \AND
  \bauthor{\bsnm{Lupparelli},~\bfnm{Monica}\binits{M.}}
(\byear{2011}).
\btitle{Chain graph models of multivariate regression type for categorical
  data}.
\bjournal{Bernoulli}
\bvolume{17}
\bpages{827--844}.
\end{barticle}
\endbibitem

\bibitem[\protect\citeauthoryear{Mat\'u\v{s}}{1992}]{matus:92}
\begin{barticle}[author]
\bauthor{\bsnm{Mat\'u\v{s}},~\bfnm{Franti\v{s}ek}\binits{F.}}
(\byear{1992}).
\btitle{On Equivalence of {M}arkov Properties over Undirected Graphs}.
\bjournal{J.\ Appl.\ Prob.}
\bvolume{29}
\bpages{745--749}.
\end{barticle}
\endbibitem

\bibitem[\protect\citeauthoryear{Pe\~{n}a}{2014}]{pen14}
\begin{barticle}[author]
\bauthor{\bsnm{Pe\~{n}a},~\bfnm{Jose~M.}\binits{J.~M.}}
(\byear{2014}).
\btitle{Marginal {AMP} chain graphs}.
\bjournal{Internat.\ J.\ Approx.\ Reason.}
\bvolume{55}
\bpages{1185--1206}.
\end{barticle}
\endbibitem

\bibitem[\protect\citeauthoryear{Pe\~{n}a}{2016}]{pen16}
\begin{binproceedings}[author]
\bauthor{\bsnm{Pe\~{n}a},~\bfnm{Jose~M.}\binits{J.~M.}}
(\byear{2016}).
\btitle{Alternative {M}arkov and causal properties for acyclic directed mixed
  graphs}.
In \bbooktitle{Proceedings of the Thirty-Second Conference on Uncertainty in
  Artificial Intelligence, {UAI} 2016, June 25-29, 2016, Jersey City, NJ,
  {USA}}.
\end{binproceedings}
\endbibitem

\bibitem[\protect\citeauthoryear{Pearl}{1988}]{pea88}
\begin{bbook}[author]
\bauthor{\bsnm{Pearl},~\bfnm{J.}\binits{J.}}
(\byear{1988}).
\btitle{Probabilistic Reasoning in Intelligent Systems : networks of plausible
  inference}.
\bpublisher{Morgan Kaufmann Publishers}, \baddress{San Mateo, CA, USA}.
\end{bbook}
\endbibitem

\bibitem[\protect\citeauthoryear{Richardson}{2003}]{ric03}
\begin{barticle}[author]
\bauthor{\bsnm{Richardson},~\bfnm{T.}\binits{T.}}
(\byear{2003}).
\btitle{Markov properties for acyclic directed mixed graphs}.
\bjournal{Scand.\ J.\ Stat.}
\bvolume{30}
\bpages{145-157}.
\end{barticle}
\endbibitem

\bibitem[\protect\citeauthoryear{Richardson and Spirtes}{2002}]{ric02}
\begin{barticle}[author]
\bauthor{\bsnm{Richardson},~\bfnm{T.~S.}\binits{T.~S.}} \AND
  \bauthor{\bsnm{Spirtes},~\bfnm{P.}\binits{P.}}
(\byear{2002}).
\btitle{Ancestral graph {M}arkov models}.
\bjournal{Ann. Statist.}
\bvolume{30}
\bpages{962--1030}.
\end{barticle}
\endbibitem

\bibitem[\protect\citeauthoryear{Sadeghi}{2013}]{sad13}
\begin{barticle}[author]
\bauthor{\bsnm{Sadeghi},~\bfnm{Kayvan}\binits{K.}}
(\byear{2013}).
\btitle{Stable mixed graphs}.
\bjournal{Bernoulli}
\bvolume{19}
\bpages{2330-2358}.
\end{barticle}
\endbibitem

\bibitem[\protect\citeauthoryear{Sadeghi}{2016}]{sad16}
\begin{barticle}[author]
\bauthor{\bsnm{Sadeghi},~\bfnm{Kayvan}\binits{K.}}
(\byear{2016}).
\btitle{Marginalization and conditioning for LWF chain graphs}.
\bjournal{Ann. Statist.}
\bvolume{44}
\bpages{1792--1816}.
\end{barticle}
\endbibitem

\bibitem[\protect\citeauthoryear{Sadeghi and Lauritzen}{2014}]{sadl14}
\begin{barticle}[author]
\bauthor{\bsnm{Sadeghi},~\bfnm{Kayvan}\binits{K.}} \AND
  \bauthor{\bsnm{Lauritzen},~\bfnm{Steffen~L.}\binits{S.~L.}}
(\byear{2014}).
\btitle{Markov properties for mixed graphs}.
\bjournal{Bernoulli,}
\bvolume{20}
\bpages{676--696}.
\end{barticle}
\endbibitem

\bibitem[\protect\citeauthoryear{Sadeghi and Wermuth}{2016}]{sadw15}
\begin{barticle}[author]
\bauthor{\bsnm{Sadeghi},~\bfnm{Kayvan}\binits{K.}} \AND
  \bauthor{\bsnm{Wermuth},~\bfnm{Nanny}\binits{N.}}
(\byear{2016}).
\btitle{Pairwise Markov properties for regression graphs}.
\bjournal{Stat}
\bvolume{5}
\bpages{286--294}.
\end{barticle}
\endbibitem

\bibitem[\protect\citeauthoryear{Spirtes, Glymour and Scheines}{2000}]{Spi00}
\begin{bbook}[author]
\bauthor{\bsnm{Spirtes},~\bfnm{P.}\binits{P.}},
  \bauthor{\bsnm{Glymour},~\bfnm{C.}\binits{C.}} \AND
  \bauthor{\bsnm{Scheines},~\bfnm{R.}\binits{R.}}
(\byear{2000}).
\btitle{Causation, Prediction, and Search},
\bedition{2nd} ed.
\bpublisher{MIT press}.
\end{bbook}
\endbibitem

\bibitem[\protect\citeauthoryear{Spirtes, Richardson and Meek}{1997}]{spi97}
\begin{btechreport}[author]
\bauthor{\bsnm{Spirtes},~\bfnm{P.}\binits{P.}},
  \bauthor{\bsnm{Richardson},~\bfnm{T.}\binits{T.}} \AND
  \bauthor{\bsnm{Meek},~\bfnm{C.}\binits{C.}}
(\byear{1997}).
\btitle{The dimensionality of mixed ancestral graphs}
\btype{Technical Report} No. \bnumber{CMU-PHIL-83},
\bpublisher{Philosophy Department, CMU}.
\end{btechreport}
\endbibitem

\bibitem[\protect\citeauthoryear{Studen{\'y}}{1989}]{stu89}
\begin{barticle}[author]
\bauthor{\bsnm{Studen{\'y}},~\bfnm{M.}\binits{M.}}
(\byear{1989}).
\btitle{Multiinformation and the problem of characterization of conditional
  independence relations}.
\bjournal{Problems of Control and Information Theory}
\bvolume{18}
\bpages{3--16}.
\end{barticle}
\endbibitem

\bibitem[\protect\citeauthoryear{Studen{\'y}}{1998}]{stu98}
\begin{binproceedings}[author]
\bauthor{\bsnm{Studen{\'y}},~\bfnm{Milan}\binits{M.}}
(\byear{1998}).
\btitle{Bayesian networks from the point of view of chain graphs}.
In \bbooktitle{Proceedings of the Fourteenth Conference on Uncertainty in
  Artificial Intelligence}
\bpages{496--503}.
\bpublisher{Morgan Kaufmann}, \baddress{San Francisco, CA}.
\end{binproceedings}
\endbibitem

\bibitem[\protect\citeauthoryear{Studen{\'y}}{2005}]{stu05}
\begin{bbook}[author]
\bauthor{\bsnm{Studen{\'y}},~\bfnm{M.}\binits{M.}}
(\byear{2005}).
\btitle{Probabilistic Conditional Independence Structures}.
\bpublisher{Springer-Verlag}, \baddress{London, United Kingdom}.
\end{bbook}
\endbibitem

\bibitem[\protect\citeauthoryear{Studen{\'y} and Bouckaert}{1998}]{stub98}
\begin{barticle}[author]
\bauthor{\bsnm{Studen{\'y}},~\bfnm{M.}\binits{M.}} \AND
  \bauthor{\bsnm{Bouckaert},~\bfnm{R.~R.}\binits{R.~R.}}
(\byear{1998}).
\btitle{On chain graph models for description of conditional independence
  structures}.
\bjournal{Ann. Statist.}
\bvolume{26}
\bpages{1434--1495}.
\end{barticle}
\endbibitem

\bibitem[\protect\citeauthoryear{Wermuth}{2011}]{wer11}
\begin{barticle}[author]
\bauthor{\bsnm{Wermuth},~\bfnm{N.}\binits{N.}}
(\byear{2011}).
\btitle{Probability distributions with summary graph structure}.
\bjournal{Bernoulli}
\bvolume{17}
\bpages{845--879}.
\end{barticle}
\endbibitem

\bibitem[\protect\citeauthoryear{Wermuth, Cox and Pearl}{1994}]{wer94}
\begin{btechreport}[author]
\bauthor{\bsnm{Wermuth},~\bfnm{N.}\binits{N.}},
  \bauthor{\bsnm{Cox},~\bfnm{D.~R.}\binits{D.~R.}} \AND
  \bauthor{\bsnm{Pearl},~\bfnm{J.}\binits{J.}}
(\byear{1994}).
\btitle{Explanation for multivariate structures derived from univariate
  recursive regressions}
\btype{Technical Report} No. \bnumber{94(1)},
\bpublisher{Univ. Mainz},
\baddress{Germany}.
\end{btechreport}
\endbibitem

\bibitem[\protect\citeauthoryear{Wermuth and Cox}{1998}]{wer98}
\begin{barticle}[author]
\bauthor{\bsnm{Wermuth},~\bfnm{N.}\binits{N.}} \AND
  \bauthor{\bsnm{Cox},~\bfnm{D.~R.}\binits{D.~R.}}
(\byear{1998}).
\btitle{On association models defined over independence graphs}.
\bjournal{Bernoulli}
\bvolume{4}
\bpages{477--495}.
\end{barticle}
\endbibitem

\bibitem[\protect\citeauthoryear{Wermuth and Lauritzen}{1983}]{wer83}
\begin{barticle}[author]
\bauthor{\bsnm{Wermuth},~\bfnm{N.}\binits{N.}} \AND
  \bauthor{\bsnm{Lauritzen},~\bfnm{S.~L.}\binits{S.~L.}}
(\byear{1983}).
\btitle{Graphical and recursive models for contingency tables}.
\bjournal{Biometrika}
\bvolume{70}
\bpages{537--552}.
\end{barticle}
\endbibitem

\bibitem[\protect\citeauthoryear{Wermuth, Marchetti and Cox}{2009}]{wer09}
\begin{barticle}[author]
\bauthor{\bsnm{Wermuth},~\bfnm{N.}\binits{N.}},
  \bauthor{\bsnm{Marchetti},~\bfnm{G.~M.}\binits{G.~M.}} \AND
  \bauthor{\bsnm{Cox},~\bfnm{D.~R.}\binits{D.~R.}}
(\byear{2009}).
\btitle{Triangular systems for symmetric binary variables}.
\bjournal{Electron. J. Stat.}
\bvolume{3}
\bpages{932--955}.
\end{barticle}
\endbibitem

\bibitem[\protect\citeauthoryear{Wermuth and Sadeghi}{2012}]{wers11}
\begin{barticle}[author]
\bauthor{\bsnm{Wermuth},~\bfnm{Nanny}\binits{N.}} \AND
  \bauthor{\bsnm{Sadeghi},~\bfnm{Kayvan}\binits{K.}}
(\byear{2012}).
\btitle{Sequences of regressions and their independences}.
\bjournal{TEST}
\bvolume{21}
\bpages{215--252 and 274--279}.
\end{barticle}
\endbibitem

\bibitem[\protect\citeauthoryear{Wright}{1921}]{wright:21}
\begin{barticle}[author]
\bauthor{\bsnm{Wright},~\bfnm{S.}\binits{S.}}
(\byear{1921}).
\btitle{Correlation and Causation}.
\bjournal{J. Agricultural Res.}
\bvolume{20}
\bpages{557--585}.
\end{barticle}
\endbibitem

\end{thebibliography}
